\documentclass[journal,twoside,web]{ieeecolor}
\usepackage{generic}
\usepackage{cite}
\usepackage{amsmath}

\usepackage{amsthm}
\usepackage{algorithm}
\usepackage{tikz}
\usetikzlibrary{tikzmark}

\usepackage{graphicx}
\usepackage{textcomp}
\usepackage{mathtools}
\usepackage{multirow}
\usepackage{capt-of}

\let\labelindent\relax
\usepackage{enumitem}
\usepackage{amssymb}
\usepackage{subcaption}
\usepackage{booktabs}
\usepackage{placeins}
\usepackage{blindtext}
\usepackage{nicefrac}
\usepackage[noend]{algpseudocode}
\usepackage[hidelinks]{hyperref}

\newtheorem{theorem}{Theorem}
\newtheorem{lemma}{Lemma}
\newtheorem{corollary}{Corollary}[theorem]
\newtheorem*{remark}{Remark}

\newtheorem{assumption}{Assumption}

\usepackage[normalem]{ulem}
\usepackage{soul}

\def\BibTeX{{\rm B\kern-.05em{\sc i\kern-.025em b}\kern-.08em T\kern-.1667em\lower.7ex\hbox{E}\kern-.125emX}}

\DeclareMathOperator*{\subjto}{subj.\,to}

\DeclareMathOperator*{\argmin}{arg\,min}

\DeclareRobustCommand{\rchi}{{\mathpalette\irchi\relax}}
\newcommand{\irchi}[2]{\raisebox{\depth}{$#1\chi$}}

\begin{document}

\thispagestyle{empty}
\onecolumn
\hspace{0em} \vfill
\noindent \textbf{This work has been submitted to the IEEE for possible publication. Copyright may be transferred without notice, after which this version may no longer be accessible.}

\vspace{1em}
\noindent First submitted to the IEEE Transactions on Automatic Control on 8/17/2025. \\
\noindent Date of latest revision 7/29/2026.
\vfill \hspace{0em}
\twocolumn
\newpage
\twocolumn
\clearpage
\pagenumbering{arabic}

\title{Observed Control---Linearly Scalable Nonlinear Model Predictive Control with Adaptive Horizons}

\author{Eugene T. Hamzezadeh and Andrew J. Petruska
    \thanks{\\  ~} 
    \thanks{This work was supported in part by the ARL TBAM-CRP [W911NF-22-2-0235]. \textit{(Corresponding author: Eugene T. Hamzezadeh.)}}%
    \thanks{Eugene T. Hamzezadeh \textit{(e-mail: ehamzeza@mines.edu)} and Andrew J. Petruska are with Colorado School of Mines, Golden, CO 80401, USA }
    \thanks{Andrew J. Petruska is also with Lunar Outpost, Golden, CO, 80401, USA \textit{(e-mail: andrew.petruska@lunaroutpost.com)} }
    \thanks{Digital Object Identifier [Insert Here]}
}

\maketitle

\begin{abstract}
    This work highlights the duality between state estimation and model predictive control. A model predictive controller, \textit{observed control}, is presented that uses this duality to efficiently compute control actions with linear prediction horizon length scalability. The proposed algorithms provide exceptional computational efficiency, adaptive time horizon lengths, and early optimization termination criteria. The use of Kalman smoothers as the backend optimization framework provides for a familiar implementation supported by strong theoretical guarantees. Additionally, a formulation is presented that separates linear model predictive control into purely reactive and anticipatory components, enabling \textit{any-time any-horizon observed control} while ensuring controller stability for short time horizons. Finally, the method is extended to nonlinear systems and non-quadratic cost functions to obtain locally-optimal control of complex systems while maintaining linear prediction horizon scalability and adaptive-horizon capabilities.
\end{abstract}

\begin{IEEEkeywords}
    Observed Control, Nonlinear Model Predictive Control, Estimation Control Duality, Adaptive Time Horizons, Kalman Filtering, Receding Horizon Control
\end{IEEEkeywords}

\section{Introduction}

Model predictive control (MPC) is a powerful control strategy that solves finite-horizon optimal control problems (OCPs) online, enabling flexible controllers that handle nonlinearities, multi-input multi-output systems, and generally, constraints. MPC's breadth of deployment, however, remains limited by the computational burden of repeatedly solving OCPs. This challenge is particularly acute for fast or nonlinear systems, where the controller's performance is tightly coupled to both the efficiency of the solution approach and the prediction horizon length. This work presents \textit{Observed Control}, an MPC technique that leverages duality with Kalman-style observers to efficiently solve nonlinear OCPs while dynamically determining the required prediction horizon length. A brief overview of the dominant OCP solution methodologies, horizon selection techniques, and duality are provided next.

\textit{Direct methods} for MPC discretize the OCP into a nonlinear program that techniques like sequential quadratic programming (SQP) solve, offering broad modeling flexibility and support for general constraints \cite{book_model_predictive_control_theory_computation_and_design}, but at a higher computational cost. \textit{Continuation methods} instead hot-start the underlying solver with the last solution and perform only one Newton-type iteration per update \cite{li_1989_multistep_newton_type_control_strategy, ohtsuka_2004_gmres_method_for_nl_rhc}, leveraging temporal coherence to reduce computational burden. Among these, real-time iteration (RTI) schemes are highly effective \cite{diehl_2001_real_time_optimization_for_large_scale_nonlinear_processes_phd, diehl_2005_rti_scheme_for_nonlinear_opt_feedback_ctrl, zavala_2009_advanced_step_controller}, with highly optimized tools for generating problem-specific controllers \cite{houska_2011_autogenerated_rti_nonlinear_mpc_microsecond, houska_2011_acado_toolkit, manual_ACADO_toolkit, verschueren_2018_acados_sofware_paper}. However, like direct methods, RTI methods typically remain quadratically scalable with respect to the prediction horizon length, with one notable exception \cite{frison_2020_hpipm_linear_qp}.

In contrast to direct optimization, \textit{indirect methods} solve OCPs by exploiting dynamic programming structure, computing solutions via forward-backward recursions. Iterative linear quadratic regulation (iLQR) \cite{li_2004_iterated_lqr, tassa_2012_ilqg} and differential dynamic programming (DDP) \cite{book_jacobson_1970_differential_dynamic_programming, tassa_2014_control_limited_ddp} are canonical methods of this type. DDP uses second-order dynamics information, while iLQR relies on first-order linearizations. These approaches combine a forward roll out of the system dynamics with a backward dynamic programming recursion, finding solutions with linear computational complexity in the horizon per iteration. However, their performance can be sensitive to initialization resulting in convergence to poor local minima.

Alternatively, \textit{sampling methods} use Monte Carlo techniques to explore the state and control spaces of OCPs without requiring gradients, enabling application to nonlinear systems and general objectives. Model predictive path integral (MPPI) control \cite{williams_2016_mppi_aggresive_driving} and variants \cite{williams_2018_robust_sampling_based_mpc, gandhi_2021_robust_mppi_analysis_perf_guarantees, wang_2021_mppi_tsallis, balci_2022_constraints_covariance_tube_mppi}
are widely used; to operate, they sample control sequences, evaluate resulting trajectories, and update the control horizon with a path-integral \cite{kappen_2005_linear_theory_for_control_of_nonlinear_stochastic_systems__path_integral}. However, a key limitation is the exponential curse of dimensionality in obtaining adequate samples, which is exasperated for unstable systems making MPPI generally impractical for systems without significant parallel computing, i.e, a GPU.

As with the solution method, the prediction horizon length $N$ is critical to computational effort and closed-loop performance. Terminal costs and constraints can shorten horizons and aid stability and optimality proofs \cite{grune_2008_on_infhoriz_perf_rhc, keerthi_1988_optimal_infinite_horizon_laws, hu_2002_toward_infinite_horizon_opt_in_nlmpc, grimm_2005_mpc_for_want_of_lyapunov_all_not_lost}, but designing terminal constraints for nonstationary references is often intractable, and terminal costs may distort the cost function with respect to the infinite-time original. Stability guarantees have been established using sufficiently large $N$ with \cite{primbs_2000_feas_and_stability_of_finite_rhc} and without \cite{jadbabaie_2005_on_rhc_stability_with_general_term_cost, boccia_2014_stability_and_feasibility_without_term_constr} terminal costs. Similarly, suboptimality bounds have been derived both with \cite{primbs_2000_feas_and_stability_of_finite_rhc, grune_2008_on_infhoriz_perf_rhc} and without \cite{shamma_1997_linear_nonquadratic_optimal_control, grune_2008_on_infhoriz_perf_rhc} terminal costs for linear and nonlinear systems. However, these approaches select a fixed $N$ offline, reusing it for each optimization resulting in bounds that are often conservative for stabilized systems and do not adapt $N$ to the current state, reference trajectory, or plant linearization.

The presented approach leverages estimation-control duality to derive the control-space dual of recursive smoothers. While prior work has extended the classical Kalman-Bucy duality \cite{kalman_1960_filter, kalmanbucy_1961_new_results_filtering_theory} through information filters \cite{torodov_2008_duality_infofilter}, backward stochastic differential equations \cite{kim_2024_duality_1_observability, kim_2024_duality_2_optctrl}, and variational principles \cite{mitter_2003_variational_nonlinear_est}, those efforts primarily deepen nonlinear duality theory. In contrast, the presented duality relationship manifests directly at the algorithmic level: Kalman-style machinery forms the core of the presented MPC algorithms and familiar estimation constructs acquire natural interpretations in predictive control. 
The development in this paper is restricted to unconstrained MPC, with extension of the algorithmic duality to general state and input constraints left for future work.

The primary contribution of this work is the development of \textit{Observed Control}, a model predictive control framework derived from recursive state estimation. Using this duality, four MPC algorithms based on Kalman-style recursions are presented that scale linearly with prediction horizon length, yielding computational complexity comparable to state-of-the-art optimization based approaches and scalability in-line with iLQR and DDP. The framework naturally produces nonlinear and derivative-free controllers by leveraging estimators such as the extended and unscented Kalman filters (EKF, UKF) \cite{book_grewal_2015_kalman_filtering, simon_2006_optimal_state_estimation,julier_1997_ukf_a_new_kf_extension}, respectively, and locally-optimal controllers with iterated variants \cite{book_sarkka_2023_bayesian_filtering_and_smoothing}. A second contribution is a formulation that computes controls using a single forward-in-time recursion, and termination metrics that enable adaptive horizon selection. A third contribution is the demonstration that linear-quadratic MPC policies are superpositions of purely reactive and anticipatory components, with the later computable offline, enabling any-time any-horizon operation. Finally, the proposed methods are evaluated on linear and nonlinear systems with quadratic and non-quadratic objectives. Although an early version of observed control was presented by the authors for magnetic devices \cite{pratt_2023_observed_control}; this work provides the first formal presentation of the approach together with its theoretical analysis, adaptive-horizon formulation, and any-time any-horizon capabilities.

The rest of the paper is structured as follows. Following preliminaries, Section~\ref{section:observed_control} develops the theory behind \textit{Observed Control} with implications for MPC in general. The treatment begins with a basic forward-backward recursion algorithm in Section~\ref{subsec:oc_naive_smoother}, which is simplified to a forward-only approach in Section~\ref{subsec:oc_foward_only}, and then to a more numerically stable EKF-specific variant in Section~\ref{subsec:oc_numerical_improvements}. The forward-only nature of the algorithms enables adaptive prediction horizons based on quantitative metrics for reference tracking performance and infinite-horizon equivalence that are described in Section~\ref{subsec:oc_early_term_crit}. Restricting the analysis to linear time-invariant systems with quadratic costs, Section~\ref{subsec:oc_anytime} proves that an optimal linear-quadratic MPC policy is separable into purely reactive and anticipatory components, ensuring closed-loop stability even for minimal one-step horizons. The efficiency and versatility of the presented algorithms are explored in case studies on the predictive control of linear and nonlinear systems with quadratic and non-quadratic objective functions in Section \ref{section:case_studies_and_disc}. For readers interested in implementation, the authors recommend starting with Algorithm~\ref{alg:oc_foward_only} (Section~\ref{subsec:oc_foward_only}), as it is compatible with both EKF and UKF computations.

\section{Preliminaries} \label{section:preliminaries}

Capital letters $A$ denote matrices; lower case $x$, vectors. Hats $\hat{x}$ indicate estimates; dots $\dot{x}$, $\ddot{x}$ time derivatives. The identity and zero matrices are $\mathbb{I}$ and $\mathbb{O}$. Subscripts $x_{k(-)}$, $x_{k(+)}$, $x_{k[s]}$ denote prior, posterior, and smoothed states at time step $k$ respectively. The shorthand $\|x\|_A^2$ is $x^\top A x$, and $A^\dagger$ a lease-squares pseudo-inverse, e.g, $(A^\top A)^{-1} A^\top$. The normal distribution with mean $\mu$ and covariance $\Sigma$ is $\mathcal{N}(\mu, \Sigma)$. All norms are $\ell^2$, written $\|A\|$ without subscripts. Bracketed $u_0^{[N]}$ denotes control $u_0$ from a prediction horizon of length $N$.
Products $\prod_{i=a}^{b}$ whose ranges satisfy $b<a$ are the identity $\mathbb{I}$.

\subsection{Model Predictive Control} \label{section:preliminaries_mpc}

This works considers MPC strategies that repeatedly solve finite-horizon discrete optimal control problems of the form 
\begin{subequations}
    \label{eqn:disc_time_opt_ctrl_formulation}
    \begin{align}
        \label{eqn:disc_time_opt_ctrl_formulation:obj_fn}
        \argmin_{u_0, \cdots, u_{N-1}}~~& \sum_{k=0}^{N} \left( J_k(x_k, u_k) + \|u_k - u_{k-1}\|_{\tilde{R}_p}^2 \right)\\
        \label{eqn:disc_time_opt_ctrl_formulation:dynamics}
        \subjto~~& x_k = f(x_{k-1}, u_{k-1})\\
        \label{eqn:disc_time_opt_ctrl_formulation:init}
        & x_0 = x(0).
    \end{align}
\end{subequations} 
The objective function \eqref{eqn:disc_time_opt_ctrl_formulation:obj_fn} is assumed to consist of the sum of generic stage terms $J_k(x_k, u_k)$ and quadratic control change penalties. It is also required to admit at least one stationary minimizer. If a terminal cost is desired, it is assumed to be implemented in $J_N(x_N, u_N)$. However, this work does not impose a terminal cost or constraint explicitly, instead the prediction horizon length $N$ is selected to sufficiently approximate the infinite-horizon problem. Throughout this work, the plant dynamics \eqref{eqn:disc_time_opt_ctrl_formulation:dynamics} are treated as the only hard constraint; general state and input constraints may be approximated via standard penalty or barrier methods.

\subsection{State Estimation} \label{section:preliminaries_state_est}

State estimators, i.e., observers, determine the instantaneous state $\hat{x}_k$ or sequence of states $(\hat{x}_0, \cdots, \hat{x}_N)$ that best explain the trajectory of a system given a noisy measurements $(y_0, \cdots, y_N)$. 
For clarity, the presentation is temporarily restricted to linear systems and measurements; Section \ref{sec:nonlinear_systems_extensions} provides a nonlinear treatment. Consider the linear system
\begin{subequations} 
    \begin{align}
        \label{linear_meas_restrictions}
        x_k &= Ax_{k-1} + Bu_{k-1} + w & w \sim \mathcal{N}(0, Q)\\
        y_k &= Cx_k + v                & v \sim \mathcal{N}(0, R)
    \end{align}
\end{subequations}
with states $x \in \mathbb{R}^n$, controls $u_k\in \mathbb{R}^m$, measurements  $y \in \mathbb{R}^\nu$, and a discretization period of $\Delta t = t_k-t_{k-1}$. The process noise $w$ and measurement noise $v$ terms are zero-mean normal distributions with covariances
$Q \in \mathbb{R}^{n \times n}$ and
$R \in \mathbb{R}^{\nu \times \nu}$.

The following algorithmic development centers around the Kalman filter and Rauch-Tung-Striebel (RTS) smoother; a brief review of the well-known equations is therefore provided \cite{kalman_1960_filter, rauch_1965_rts_smoother}. The Kalman filter's mean and covariance prediction equations are given by
\begin{subequations} \label{eq: forward pass}
    \begin{align}
        \label{kf:dynamics}
        \hat{x}_{k(-)} &= A \hat{x}_{k-1(+)} + B u_{k-1}\\
        \label{kf:prior_propagation}
        P_{k(-)} &= A P_{k-1(+)} A^\top + Q.
    \end{align}
When a new measurement $y_k$ arrives, the anticipated measurement residual $r_k$ is given by 
    \begin{equation}
        \label{kf:anticipated_residual}
        r_k = y_k - C\hat{x}_{k(-)}.
    \end{equation}
The posterior prediction and covariance are produced with
    \begin{align}
        \label{kf:posterior_state_update}
        \hat{x}_{k(+)} &= \hat{x}_{k(-)} + K_k r_k\\
        \label{kf:posterior_cov_update}
        P_{k(+)} &= (\mathbb{I} - K_k C) P_{k(-)}
    \end{align}
where $K_k$ is the Kalman gain given by
    \begin{align}
        \label{kf:gain}
        K_k &= P_{k(-)}C^\top (CP_{k(-)}C^\top + R)^{-1}.
    \end{align}
\end{subequations}
Smoothers improve the estimate for $\hat{x}_k$ by incorporating future measurements $(y_{k+1}, \cdots, y_N)$. The RTS smoother is one of the most commonly used Kalman smoothers, which unlike previous smoothers, is formulated to only use information from the forward filter pass. Also expressed recursively, operating backwards in time, the RTS smoother equations are
\begin{subequations}
    \label{rts_equations}
    \begin{align}
        \label{eq:kalman_smoother_state}
        \hat{x}_{k[s]} &= \hat{x}_{k(+)} + L_k (\hat{x}_{k+1[s]} - \hat{x}_{k+1(-)}) \\
        \label{eq:kalman_smoother_covariance}
        P_{k[s]} &= P_{k(+)} + L_k (P_{k+1[s]} - P_{k+1(-)}) L_k^\top \\
        \label{eq:kalman_smoother_gain}
        L_k &= P_{k(+)} A_k^\top P_{k+1(-)}^{-1}
    \end{align}
\end{subequations}
where $L_k$ is the smoother gain. The smoother is initialized by setting $\hat{x}_{N[s]}$ to the last posterior $\hat{x}_{N(+)}$. Then, to produce the set of smoothed state estimates $( \hat{x}_{0[s]}, \cdots, \hat{x}_{N[s]} )$, equations (\ref{eq:kalman_smoother_state}-\ref{eq:kalman_smoother_gain}) are iterated backwards to $k = 0$. 

\section{Observed Control} \label{section:observed_control}

The central idea behind observed control is to reformulate predictive control problems as state estimation tasks. The framework uses a forward filter pass to predict future states and construct synthetic measurements that encode the desired objective function. A subsequent smoother pass then estimates the controls that best explain the measurements. From a maximum-likelihood perspective, the smoother finds the controls that are most likely to produce the desired observations, i.e., the those corresponding to the objective function. In this reformulation, dynamics constraints are encoded in the filter's process model while the objective function is encoded in the measurement residuals. Under an appropriate mapping, the resulting estimation problem is equivalent to the original predictive control problem, with the objective function appearing as the negative log-likelihood minimized by the smoother. Notably, this estimation-based perspective on predictive control enables several computational and theoretical advantages.

The following development maps the linear-quadratic model predictive control problem into an equivalent smoother formulation, setting up the augmentation concepts of the observed control framework. At this point, the analysis remains restricted to linear time-invariant systems of the form
\begin{subequations}
    \begin{align}
        x_{k} &= Ax_{k-1} + Bu_{k-1}
    \end{align}
\end{subequations}
with $n$ states and $m$ controls. To use the observed control framework to \textit{estimate controls}, the number of augmented states $\eta$ is defined as $\eta = n + m$. The new augmented system state $\rchi_k$ and state transition matrix $\Phi$ are defined as
\begin{align}
    \label{augmented_system}
    \rchi_k = \begin{bmatrix}
        x_k\\
        u_k
    \end{bmatrix}&&
    \Phi = \begin{bmatrix}
        A & B\\
        \mathbb{O} & \mathbb{I}_{m}
    \end{bmatrix}
\end{align}
where $\rchi \in \mathbb{R}^\eta$ and $\Phi \in \mathbb{R}^{\eta \times \eta}$. The augmented model with $\mu$ augmented measurements used for filtering and smoothing is
\begin{subequations}
    \label{eqn:linear_systems_augmented_model}
    \begin{align}
        \hat{\rchi}_k &= \Phi \hat{\rchi}_{k-1} + \grave{w} & \grave{w} \sim \mathcal{N}(0, \grave{Q})\\
        \label{eqn:linear_systems_augmented_model:just_the_measurement}
        z_k &= H \hat{\rchi}_k + \grave{v}            & \grave{v} \sim \mathcal{N}(0, \grave{R})
    \end{align}
\end{subequations}
with $z_k \in \mathbb{R}^\mu$, $H\in \mathbb{R}^{\mu \times \eta}$, $\grave{Q} \in \mathbb{R}^{\eta \times \eta}$, and $\grave{R} \in \mathbb{R}^{\mu \times \mu}$. 
The construction of the measurement model \eqref{eqn:linear_systems_augmented_model:just_the_measurement} is covered next.
To avoid notation clutter, the circumflex $(\hat{~})$ character will henceforth be omitted on the augmented state estimate $\rchi$.

\subsection{Linear-Quadratic MPC Equivalence} \label{subsec:oc_reduction_to_lqr}
Given the well-known duality between the Kalman filter and LQR, it only seems natural to explain the mapping of the objective and constraints within the observed control framework using a linear-quadratic MPC formulation to start. This will inform the construction of the filter's covariance matrices and measurements. Consider reference tracking formulation
\begin{subequations}
    \label{eqn:lqr_forulation}
    \begin{alignat}{1}
        \label{eqn:lqr_minimization}
        \argmin_{(u_0, \cdots, u_N)} ~~& \sum_{k=0}^{N} \Bigl( \|x_{r, k} -x_k\|^2_{Q_{p}} + \|u_{k}\|^2_{R_{p}} \\
        \notag
            & \qquad  +\:2x_k^T M_{p} u_k + \|u_k-u_{k-1}\|^2_{\tilde{R}_{p}} \Bigr) \\
        \label{eqn:lqr_dynamics_constraint}
        \subjto ~~& x_k=Ax_{k-1}+Bu_{k-1} \\
        \label{eqn:lqr_init_constraint}
        & x_0 = x(0)
    \end{alignat}
\end{subequations}
where $x_{r,k}$ is the desired reference at index $k$, $x(0)$ is the current system state, $u_{-1}$ is the previously applied control, and $Q_{p}$, $R_{p}$, $M_{p}$, and $\tilde{R}_{p}$ are weighting matrices.

\begin{assumption}
    \label{assumpion_reasonable_weights_detectable_sys}
    The objective function weighting matrices satisfy $Q_{p} \succ 0$, $R_{p} \succ 0$, and $\tilde{R}_{p} \succ 0$. The cross-term weighting matrix $M_{p}$ satisfies $Q_p - M_p R_p^{-1} M_p^\top \succeq 0$. The system $x_k=Ax_{k-1}+Bu_{k-1}$ is stabilizable.
\end{assumption}

\begin{theorem} \label{thm:oc_lqr_equilvalence}
    Let Assumption \ref{assumpion_reasonable_weights_detectable_sys} hold. Then an augmented RTS smoother defined by the system \eqref{eqn:linear_systems_augmented_model}, measurements \eqref{eqn:oc_lqr_synth_measurement}, and covariances \eqref{eqn:oc_lqr_qr_mapping} produces the unique optimizer of the finite-horizon problem \eqref{eqn:lqr_forulation}. Furthermore as $N \rightarrow \infty$, the smoother solution converges to the infinite-horizon solution of \eqref{eqn:lqr_forulation}.
\end{theorem}

\begin{proof}
    The RTS smoother optimally minimizes the negative log-likelihood of the maximum a posteriori \cite{bell_1994_iterated_kalman_smoother_as_gauss_newton_method}, which for the augmented system \eqref{eqn:linear_systems_augmented_model} is
    \begin{subequations}
        \label{eqn:eks_minimization}
        \begin{alignat}{2}
            \label{eqn:eks_minimization:first_sum}
            J(\rchi_0, \cdots, \rchi_N) & =\: && \|(\rchi_0-\bar{\rchi}_0)\|_{P_{0}^{-1}}^2 \\
            \label{eqn:eks_minimization:second_sum}
            & && + \sum_{k=0}^N \|z_k - H\rchi_k\|_{R^{-1}}^2 \\
            \label{eqn:eks_minimization:third_sum}
            & && + \sum_{k=1}^{N} \|\rchi_{k} - \Phi\rchi_{k-1} \|_{Q^{-1}}^2
        \end{alignat}
    \end{subequations}
    where the filter's initialization $\bar{\rchi}_0$, synthetic measurements $z_k$ and measurement sensitivities $H$ are chosen as
    \begin{align}
        \label{eqn:oc_lqr_synth_measurement}
        \bar{\rchi}_0^\top = [ x(0)^\top ~ u_{-1}^\top ]
        &&
        z_k = \begin{bmatrix}
            x_{r,k} \\
            \mathbb{O}_{m\times1}
        \end{bmatrix}
        &&
        H=\mathbb{I}.
    \end{align}
    Choosing the smoother covariances $P_0^{-1}$, $\grave{Q}^{-1}$ and $\grave{R}^{-1}$ as 
    \begin{align*}
        P_0^{-1} = \grave{Q} ^{-1}
        &&
        \grave{Q}^{-1} = \begin{bmatrix}
            \beta \mathbb{I}_{n} & 0\\
            0 & \tilde{R}_{p}
        \end{bmatrix}
        &&
        \grave{R}^{-1} = \begin{bmatrix}
           Q_{p} & M_{p}\\
            M_{p}^\top &  R_{p}
        \end{bmatrix}
    \end{align*}
    and expanding \eqref{eqn:eks_minimization} with the augmented definition \eqref{augmented_system} yields
    \begin{align}
        \label{eqn:expanded_map_objective}
        J(&\rchi_0, \cdots, \rchi_N)  =\:  \| x_0 - x(0) \|^2_{\beta \mathbb{I}_{n}} + \| u_0 - u_{-1}\|^2_{\tilde{R}_{p}} \\
        \notag
        &~~ \begin{aligned}
            + \sum_{k=0}^N &\left( \| x_{r, k} - x_k \|_{Q_{p}}^2 + 2x_k^T M_{p} u_k+ \| u_k \|_{R_{p}}^2 \right)
        \end{aligned} \\
        \notag
        &~~ \begin{aligned}
            + \sum_{k=1}^N &\left( \| x_k-Ax_{k-1}-Bu_{k-1} \|_{\beta \mathbb{I}_{n }}^2 + \| u_k - u_{k-1} \|^2_{\tilde{R}_{p}} \right).
        \end{aligned}
    \end{align}
    Taking the limit $\beta \rightarrow \infty$, the initial state error and dynamics prediction errors are infinitely penalized, enforcing constraints \eqref{eqn:lqr_dynamics_constraint}-\eqref{eqn:lqr_init_constraint}. Additionally, the inverse $(\beta \mathbb{I}_{n})^{-1} \rightarrow \mathbb{O}_n$; therefore, the actual $\grave{Q}$ and $\grave{R}$ to be used in the smoother are given by
    \begin{align}
        \label{eqn:oc_lqr_qr_mapping}
        \grave{Q} = \begin{bmatrix}
            \mathbb{O}_{n} & \mathbb{O}\\
            \mathbb{O} & \tilde{R}_{p}^{-1}
        \end{bmatrix}
        &&
        \grave{R} = \begin{bmatrix}
           Q_{p} & M_{p}\\
            M_{p}^\top &  R_{p}
        \end{bmatrix}^{-1}
        && P_0 = \grave{Q}.
    \end{align}
    Combining the terms in the summations of \eqref{eqn:expanded_map_objective} yields
    \begin{align*}
        J(\rchi_0, \cdots, \rchi_N) = \sum_{k=0}^{N} \Bigl( &\|x_{r, k} - x_k\|^2_{Q_{p}} + 2x_k^T M_{p} u_k \\
            &+\: \|u_{k}\|^2_{R_{p}} + \|u_k-u_{k-1}\|^2_{\tilde{R}_{p}} \Bigr)
    \end{align*}
    and explicitly extracting the constraints in the limit results in
    \begin{align*}
        \argmin_{(u_0, \cdots, u_N)} ~~& J(\rchi_0, \cdots, \rchi_N) \\
        \subjto ~~& x_k=Ax_{k-1}+Bu_{k-1} \\
        & x_0 = x(0)
    \end{align*}
    which is equivalent with \eqref{eqn:lqr_forulation}. The synthetic measurement sensitivity \eqref{eqn:oc_lqr_synth_measurement} makes the pair $(\Phi, H)$ observable. Consequently, under Assumption \ref{assumpion_reasonable_weights_detectable_sys} the augmented estimation problem satisfies the conditions where the Kalman filter and RTS smoother recover the unique optimal control sequence.  \cite{kalman_1960_contribs_to_theory_of_optimal_ctrl, rauch_1965_rts_smoother}. 
\end{proof}

\begin{remark}
    If there are states and controls that should not be penalized, then a full-sized $\grave{R}^{-1}$ is rank deficient. Instead, the information that specific states and controls should be excluded from the optimization can be shifted from the weighting matrix $\grave{R}^{-1}$ to the measurement sensitivity $H_k$ so as not to measure them. Thus, $\grave{R}^{-1}$ should be reduced to only penalize the states and controls of interest, and $H_k$ and $z_k$ should be appropriately selected to ignore unpenalized combinations.
\end{remark}

\begin{remark}
    The weighting matrix $\tilde{R}_{p}$ penalizes $\Delta u$ and is recommended for improved conditioning. In practice, a $\Delta u$ penalty is desirable to avoid control chatter. Choosing $\tilde{R}_{p} \ll \grave{R}$ renders it effectively negligible. Alternatively, augmenting $\Phi$ in \eqref{augmented_system} with $\mathbb{O}$ instead of $\mathbb{I}$ removes the $\Delta u$ penalty.
\end{remark}

\subsection{Smoother-Derived Observed Control} \label{subsec:oc_naive_smoother}
Algorithm \ref{alg:oc_plain}---\textit{Observed Control}---implements an RTS smoother directly on the augmented system \eqref{eqn:linear_systems_augmented_model} with measurement and process noise defined by \eqref{eqn:oc_lqr_qr_mapping}. The algorithm initializes a Kalman filter with the current state and previously applied control and performs a forward filter pass, equations \eqref{kf:dynamics}-\eqref{kf:gain}, taking measurements as the desired state $z_k$ over a time horizon of length $N$ of the form \eqref{eqn:oc_lqr_synth_measurement}. It then performs the smoother equations \eqref{eq:kalman_smoother_state}-\eqref{eq:kalman_smoother_gain} back to time $k=0$. The desired control action is then contained in the smoothed state $\rchi_{0[s]}$.

\begin{algorithm}[!t]
    \caption{Observed Control}\label{alg:oc_plain}
    \begin{algorithmic}[1]
    \State {\textbf{initialize} $\rchi_{0(-)} \gets \begin{bmatrix}
        \hat{x}^\top & u_\text{last}^\top
    \end{bmatrix}^\top $, $P_{0(-)} \gets \grave{Q}$}
    \\
    \For {$k \in [0, \cdots, N]$} \Comment{Kalman filter pass}
        \If {$k > 0$}
            \label{alg1_line:kf_start}
            \State{$\rchi_{k(-)} \gets \Phi_{k-1} \rchi_{k-1(+)}$} \Comment{\eqref{kf:dynamics}}
            \State {$P_{k(-)} \gets \Phi_{k-1} P_{k-1(+)} \Phi_{k-1}^\top + \grave{Q}$} \Comment{\eqref{kf:prior_propagation}}
        \EndIf
        \State{$r_k \gets z_k - H_k \rchi_{k-1(-)}$} \label{alg1_line:kf_residual} \Comment{\eqref{kf:anticipated_residual}}
        \State {$K_k = P_{k(-)} H_k^\top [H_k P_{k(-)} H_k^\top + \grave{R}]^{-1}$} \Comment{\eqref{kf:gain}}
        \State {$\rchi_{k(+)} \gets \rchi_{k(-)} + K_kr_k$} \label{alg1_line:kf_posterior} \Comment{\eqref{kf:posterior_state_update}} 
        \State {$P_{k(+)} \gets [\mathbb{I} - K_k H_k] P_{k(-)}$} \label{alg1_line:kf_end} \Comment{\eqref{kf:posterior_cov_update}}
    \EndFor

    \For {$k \in [N-1, \cdots, 0]$} \Comment{RTS smoother pass}
        \State {$L_k = P_{k(+)} \Phi_k^\top P_{k+1(-)}^{-1} $} \Comment{\eqref{eq:kalman_smoother_gain}}
        \State {$\rchi_{k[s]} \gets \rchi_{k(+)} + L_k (\rchi_{k+1[s]} - \rchi_{k+1(-)})$}  \Comment{\eqref{eq:kalman_smoother_state}}
        \State {$P_{k[s]} \gets P_{k(+)} + L_k (P_{k+1[s]} - P_{k+1(-)}) L_k^\top $} \Comment{\eqref{eq:kalman_smoother_covariance}}
    \EndFor\\
    \Return $u_0 \in \rchi_{0[s]}$
    \end{algorithmic}
\end{algorithm}

\subsection{Forward Only General Observed Control} \label{subsec:oc_foward_only}
While Algorithm \ref{alg:oc_plain} is sufficient to compute the optimal predictive control update, it leaves room for computational and numerical stability improvements. In the MPC scheme, only the first control in the optimization horizon is ever applied to the system before the next optimization cycle. That is, refinements to the controls in the intermediate smoothed results $( \rchi_{1[s]} \cdots \rchi_{N[s]} )$ need not be explicitly computed. Re-examining the RTS smoother state update equations with this insight, the initial control $\rchi_{0[s]}$ can be extracted by recursively substituting equation \eqref{eq:kalman_smoother_state} with future versions of itself, yielding the relationship
\begin{align}
    \notag
    \rchi_{0[s]}
    &~\begin{aligned} 
         = \rchi_{0(+)}+L_{0}(\rchi_{1[s]}-\rchi_{1(-)})
    \end{aligned}\\
    \notag
    &~\begin{aligned}
        = \rchi_{0(+)}+L_{0}([\rchi_{1(+)}+L_{1}(\rchi_{2[s]}-\rchi_{2(-)})]-\rchi_{1(-)})
    \end{aligned}\\
    &~\begin{aligned}
        = \rchi_{0(+)} & + L_{0}(\rchi_{1(+)}-\rchi_{1(-)})\\
        & + L_{0}L_{1}(\rchi_{2(+)}-\rchi_{2(-)})\\
        & + L_{0}L_{1}L_{2}(\rchi_{3(+)}-\rchi_{3(-)}) + \{\cdots\}
    \end{aligned}
\end{align}
which highlights the individual contributions of each future measurement in the receding horizon $(k \ge 1)$ to the present-time smoothed state $(k=0)$.
In general, for $N$ measurements, the smoothed augmented state can be described by
\begin{subequations}
    \begin{align}
        \label{smoothed_state_calculation_lti}
        \rchi_{0[s]} &= \rchi_{0(+)} + \sum_{k=1}^{N} \Big( \prod_{i=0}^{k-1}L_{i} \Big) (\rchi_{k(+)}-\rchi_{k(-)}) \\
        \label{smoothed_state_calculation_lti_gamma}
        u_{0[s]} &= u_{0(+)} + \sum_{k=1}^{N} \Gamma_{k-1} (\rchi_{k(+)}-\rchi_{k(-)})
    \end{align}
\end{subequations}
where the accumulator $\Gamma_k$ can be recursively found with
\begin{align}
    \Gamma_k = \Gamma_{k-1} L_{k} &&
    \Gamma_{-1}\equiv \begin{bmatrix} \mathbb{O} &\mathbb{I}\end{bmatrix}_{m\times \eta}    .\label{gain_accumulator_eqn}
\end{align}
The matrix $\Gamma_{-1}$ can be initialized in this manner because the initial process noise $P_0$ has no uncertainty in the states $x_0 \in \rchi_0$, see \eqref{eqn:oc_lqr_qr_mapping} in Theorem \ref{thm:oc_lqr_equilvalence}.
Thus, the product $\prod_{i=0}^{k-1}L_{i}$ is identically zero for the top $n$ rows corresponding to updates in the initial states $x_0$, which need not be explicitly carried out during the computation. 

The relationship \eqref{smoothed_state_calculation_lti_gamma} removes the need for a separate smoother pass as the gain product $\Gamma_k$ can be accumulated with information only from the forward pass. A similar relationship for the smoothed initial controls covariance $P_{0[s],u}$ can be found in the same way, yielding a forward-only update on both the initial control and uncertainty at each step $k$
\begin{align}
    \label{smoothed_state_calculation_simplified}
    u_{0[s]} & \gets u_{0[s]} + \Gamma_{k-1}(\rchi_{k(+)}-\rchi_{k(-)})\\
    \label{smoother_cov_forward_simplified}
    P_{0[s],u} & \gets P_{0[s],u} + \Gamma_{k-1}(P_{k(+)} - P_{k(-)})\Gamma_{k-1}^\top.
\end{align}
Using \eqref{smoothed_state_calculation_simplified}--\eqref{smoother_cov_forward_simplified}, Algorithm \ref{alg:oc_foward_only}---\emph{Forward-Only Observed Control}---is developed, which does not spend time refining future controls $( u_1, \cdots, u_{N-1} )$ that will be discarded, and only smooths controls, not states.

The accumulator $\Gamma_k$ defines how future measurements (cost function evaluations) map to the current control action and informs the rate of convergence in the initial state uncertainty, or covariance. Understanding how this gain converges is the topic of Section \ref{subsec:oc_forward_pass_convergence} and informs adaptive time horizons and early termination metrics discussed in Section \ref{subsec:oc_early_term_crit}, but are included in Algorithm \ref{alg:oc_foward_only} for implementation completeness.

\subsection{Convergence of Forward Pass} \label{subsec:oc_forward_pass_convergence}

The control action $u_{0[s]}$ and uncertainty $P_{0[s],u}$ are dependent on the convergence of the series in \eqref{smoothed_state_calculation_simplified} and \eqref{smoother_cov_forward_simplified} which are ultimately dominated by the smoother gain product $\Gamma_k$ in \eqref{gain_accumulator_eqn}.
As $k \rightarrow \infty$, it can be shown that $\Gamma_k \rightarrow 0$, which yields useful properties for informing appropriate horizon lengths.

\begin{assumption}
    \label{assumption_kf_converges}
    The pair $(\Phi, H)$ is detectable and the pair $(\Phi, \grave{Q}^{1/2})$ is stabilizable. Consequently, the associated Kalman filter converges to a unique steady-state solution with finite covariance matrices $P_{\infty(-)}$, $P_{\infty(+)}$, and stabilizing gain $K_\infty$.
\end{assumption}

\begin{lemma} \label{lem: L seq converges} 
    Let Assumption \ref{assumption_kf_converges} hold. The sequence of smoother gains $(L_0, L_1, \cdots, L_N)$ converges to an infinite time value of 
    \begin{equation*}
        L_\infty=P_{\infty(+)} \Phi^\top P_{\infty(+)}^{-1} (\mathbb{I} - K_\infty H).
    \end{equation*}
\end{lemma}

\begin{proof}
    Under Assumption \ref{assumption_kf_converges}, the Kalman gain and covariance matrices converge to $K_\infty$, $P_{\infty(-)}$, and $P_{\infty(+)}$. The steady-state posterior covariance relation relation \eqref{kf:posterior_cov_update} implies
    \begin{equation}
        P_{\infty(-)}^{-1}=P_{\infty(+)}^{-1} (\mathbb{I} - K_\infty H)
    \end{equation}
    which can be substituted into the infinite-time version of \eqref{eq:kalman_smoother_gain}
    \begin{equation}
        L_\infty = P_{\infty(+)} \Phi^\top P_{\infty(-)}^{-1}
    \end{equation}
    yielding the desired result
    \begin{equation}
        L_\infty = P_{\infty(+)} \Phi^\top P_{\infty(+)}^{-1} (\mathbb{I} - K_\infty H)
    \end{equation}
    which completes the proof.
\end{proof}

\begin{algorithm}[!t]
    \caption{Forward-Only Observed Control}\label{alg:oc_foward_only}
    \begin{algorithmic}[1]
    \State {\textbf{initialize} $\rchi_{0(-)} \gets \begin{bmatrix}
        \hat{x}^\top & u_\text{last}^\top
    \end{bmatrix}^\top$, $u_{0[s]} \gets u_\text{last}$, }
    \State{\qquad \qquad $P_{0(-)} \gets \grave{Q}$, $P_{0[s]}\gets\grave{Q}$, $\Gamma_{-1}\gets \begin{bmatrix} \mathbb{O} & \mathbb{I} \end{bmatrix}$}
    \\
    \For {$k \in [0, \cdots, N]$} \Comment{Single Forward Pass}
        \State {\textbf{perform} Algorithm \ref{alg:oc_plain} \textbf{lines} \ref{alg1_line:kf_start}-\ref{alg1_line:kf_end}}
        \State {$\Gamma_{k-1} \gets \Gamma_{k-2} ( P_{k-1(+)} \Phi_{k-1}^{\top} P_{k(-)}^{-1} )$} \Comment{\eqref{gain_accumulator_eqn}}
        \\
        \State {$u_{0[s]} \gets u_{0[s]} + \Gamma_{k-1} (\rchi_{k(+)} - \rchi_{k(-)})$} \Comment{\eqref{smoothed_state_calculation_simplified}}
        \State {$P_{0[s],u} \gets P_{0[s],u} + \Gamma_{k-1} (P_{k(+)} - P_{k(-)}) \Gamma_{k-1}^\top $} \Comment{\eqref{smoother_cov_forward_simplified}}
        \\
        \State{$\rho_k \gets \| \Gamma_{k-1} K_k \|$} \Comment{\eqref{foward_only_rho_term}}
        \State{$\tau_k \gets  \text{tr}\left(\Gamma_{k-1} (P_{k(-)} - P_{k(+)}) \Gamma_{k-1}^\top \right) / \text{tr}(P_{0[s]})$} \Comment{\eqref{eqn:explicit_tau_compute_gammas}}
        \State{\textbf{if} $\rho_k \le \delta_1$ \textbf{and} $\tau_k \le \delta_2$ \textbf{break}} \Comment{Early Termination}
    \EndFor \\
    \Return $u_{0[s]}$
    \end{algorithmic}
\end{algorithm}

\begin{lemma}
    \label{lemma_equivalent_eigenvalues}
    Let Assumption \ref{assumption_kf_converges} hold. The matrix expressions
    \begin{align}
        \label{lemma1:L_infinity}
        L_\infty &= P_{\infty(+)} \Phi^\top P_{\infty(+)}^{-1} ( \mathbb{I} - K_\infty H) \\
        \label{lemma1:closed_loop}
        A^* &= \Phi^\top ( \mathbb{I} - H^\top K_\infty^\top )
    \end{align}
    for the infinite-time smoother gain $L_\infty$ and the given closed loop state transition matrix $A^*$, result in similar matrices.
\end{lemma}

\begin{proof}
    Two matrices $L_\infty$ and $A^*$ are similar if there exists a $V$ such that $L_\infty = V A^* V^{-1}$. Choosing $V = P_{\infty(+)}$, the product $V A^* V^{-1}$ becomes
    \begin{equation}
        V A^* V^{-1} = P_{\infty(+)} ( \Phi^\top - \Phi^\top H^\top K_\infty^\top ) P_{\infty(+)}^{-1}.
    \end{equation}
    Comparing with the expanded \eqref{lemma1:L_infinity} given by
    \begin{equation}
        L_\infty = P_{\infty(+)} \Phi^\top P_{\infty(+)}^{-1} - P_{\infty(+)} \Phi^\top P_{\infty(+)}^{-1} K_\infty H
    \end{equation}
    and explicitly setting $L_\infty = V A^* V^{-1}$ yields the equality
    \begin{equation}
        P_{\infty(+)}^{-1} K_\infty H = H^\top K_\infty^\top P_{\infty(+)}^{-1} 
    \end{equation}
    with cancellations omitted. Rearranging yields
    \begin{equation}
        K_\infty H = P_{\infty(+)} H^\top K_\infty^\top P_{\infty(+)}^{-1}.
    \end{equation}
    Substituting $P_{\infty(+)} = ( \mathbb{I} - K_\infty H) P_{\infty(-)}$, the infinite time prior-posterior relationship given by \eqref{kf:posterior_cov_update}, yields
    \begin{equation}
        K_\infty H = \left[ ( \mathbb{I} - K_\infty H ) P_{\infty(-)} \right] H^\top K_\infty^\top \left[ ( \mathbb{I} - K_\infty H ) P_{\infty(-)} \right]^{-1} 
    \end{equation}
    which, with rearranging, results in
    \begin{equation}
         K_\infty H = P_{\infty(-)} H^\top K_\infty^\top P_{\infty(-)}^{-1}.
    \end{equation}
    Substituting the infinite time Kalman gain equation \eqref{kf:gain} on the right $K_\infty = P_{\infty(-)} H^\top ( H P_{\infty(-)} H^\top + \grave{R})^{-1}$ yields
    \begin{align}
         K_\infty H &= P_{\infty(-)} \left[ P_{\infty(-)} H^\top (H P_{\infty(-)} H^\top + \grave{R})^{-1} \right]^\top P_{\infty(-)}^{-1} \nonumber \\
         &= P_{\infty(-)} (H P_{\infty(-)} H^\top + \grave{R})^{-T} H P_{\infty(-)}^\top P_{\infty(-)}^{-1}\\
         &= P_{\infty(-)} (H P_{\infty(-)} H^\top + \grave{R})^{-1} H \nonumber \\
         &= K_\infty H \nonumber
    \end{align}
    which completes the proof.
\end{proof} 

\begin{theorem} 
    \label{thm: linf_stable}
    Let Assumption \ref{assumption_kf_converges} hold. The RTS smoother gain converges to a an infinite-time value of $L_\infty$ and has all stable eigenvalues, that is $| \lambda_i | < 1$ for all $i \in [0, \eta-1]$.
\end{theorem}

\begin{proof}
    Define $e_{k(-)} = \bar{\rchi}_{k} - \rchi_{k(-)}$ and $e_{k(+)} = \bar{\rchi}_{k} - \rchi_{k(+)}$ as the Kalman filter's prior and posterior estimate errors w.r.t. the true state $\bar{\rchi}_k$. The residual $r_k$ can be expressed as
    \begin{align*}
        r_k = H\bar{\rchi}_k - H\rchi_{k(-)} = H e_{k(-)}.
    \end{align*}
    Using \eqref{kf:posterior_state_update}, the posterior update error can be written as
    \begin{align*}
        e_{k(+)} &= \bar{\rchi}_{k} -  (\rchi_{k(-)} + K_k r_k) \\
        &= (\mathbb{I} - K_k H) \Phi e_{k-1(+)}
    \end{align*}
    yielding the homogeneous error dynamics. Under Assumption \ref{assumption_kf_converges}, the error dynamics transition matrix $\Phi_e$ converges to
    \begin{equation}
        \label{eqn:error_dynamics_phi}
        \Phi_e = [\mathbb{I}-K_\infty H]\Phi
    \end{equation}
    and is guaranteed to be Schur stable. Observing that $A^*$ in \eqref{lemma1:closed_loop} is the transpose of $\Phi_e$, by Lemma \ref{lemma_equivalent_eigenvalues}, it follows that the eigenvalues $L_\infty$ are, therefore, guaranteed to be stable.
\end{proof}

\begin{lemma}
    \label{lemma_gamma_geometric_convergence}
    Let Assumption \ref{assumption_kf_converges} hold. Then, the accumulated smoother gain $\Gamma_k \rightarrow \mathbb{O}$ geometrically after a transient period.
\end{lemma}

\begin{proof}
    Consider the discrete linear time varying system
    \begin{align}
        \label{smoother_dynamics_system}
        x_k = L_k^\top x_{k-1}
    \end{align}
    which evolves as the transpose of \eqref{gain_accumulator_eqn}.
    Theorem \ref{thm: linf_stable} proves $L_\infty$, is Schur stable. Hence $L_\infty^\top$ is also Schur stable. By the discrete Lyapunov theorem, there exists a symmetric $P \succ 0$ such that
    \begin{align}
        \label{lyapunov_equation_first}
         L_\infty P L_\infty^\top - P = -Q
    \end{align}
    for a symmetric $Q \succ 0$. The Lyapunov function is defined as $V(x_k)= x_k^\top P x_k$ and strictly decreases along $x_k = L_\infty^\top x_{k-1}$. Applying the same value function to the system \eqref{smoother_dynamics_system} yields
    \begin{align}
        \label{value_function_difference}
        V(x_{k}) - V(x_{k-1}) &= x^\top_{k-1} (L_k P L_k^\top - P)x_{k-1}.
    \end{align}
    Let $L_k = L_\infty + E_k$ where $E_k$ is the convergence error term. Expanding and using \eqref{lyapunov_equation_first}, the inner term of \eqref{value_function_difference} becomes
    \begin{align}
        L_k P L_k^\top - P = L_\infty P E_k^\top + E_k P L_\infty^\top + E_k P E_k^\top - Q.
    \end{align}
    It follows from Lemma \ref{lem: L seq converges} that $E_k \rightarrow 0$ as $k \rightarrow \infty$. 
    Therefore there exists a finite $N$ where for all $k \ge N$ the relation
    \begin{align}
        L_k P L_k^\top - P \preceq -\xi Q
    \end{align}
    is satisfied for a margin $\xi \in (0, 1)$. As such the function $V(x_k)$ satisfies a strict quadratic decrease condition for $k \ge N$, and by the standard discrete-time Lyapunov exponential-stability result \cite{book_model_predictive_control_theory_computation_and_design}, the tail of the system \eqref{smoother_dynamics_system} is exponentially stable. Thus there exist constants $c>0$ and $\alpha \in(0, 1)$ such that
    \begin{align}
        \left\Vert L_k^\top L_{k-1}^\top \cdots L_{m+1} \right\Vert \le c \alpha^{k-m}, && k \ge m \ge N
    \end{align}
    is satisfied.
    Taking the transpose gives the same bound for the right-products appearing in $\Gamma_k$ and the finite transient interval $k < N$ only changes the constant $c$; Therefore, for $c_\Gamma > 0$
    \begin{align}
        \| \Gamma_k \| \le c_\Gamma \alpha ^{k-N} \| \Gamma_N\|, && k \ge N.
    \end{align}
     Hence, $\Gamma_k$ converges to $\mathbb{O}$ geometrically.
\end{proof}

\begin{corollary}  \label{cor: delta epsilon}
    Let Assumption \ref{assumption_kf_converges} hold and assume a sufficient length $N_r>N$ of bounded references $(z_0, \cdots, z_{N})$. For an arbitrarily small $\rho > 0$, there exists a prediction horizon length $N$ such that the marginal control refinement satisfies
    \begin{equation}
        \|u_0^{[N]} - u_0^{[N-1]}\| \le \rho.
    \end{equation}
\end{corollary}

\begin{proof}
    The posterior-prior difference $\rchi_{k(+)} - \rchi_{k(-)}$ can be written in terms of the Kalman gain, yielding the relation
    \begin{align}
        \label{eqn_posterior_prior_difference}
        \rchi_{k(+)} - \rchi_{k(-)} = K_k r_k.
    \end{align}
    The marginal control contribution \eqref{smoothed_state_calculation_simplified} from recursion $N$ is
    \begin{align}
        \label{marginal_control_contribution}
        u_0^{[N]} - u_0^{[N-1]} = \Gamma_{N-1} K_N r_N.
    \end{align}
    By Lemma \ref{lemma_gamma_geometric_convergence}, $\Gamma_{N-1} \rightarrow \mathbb{O}$ as as $N \rightarrow \infty$ geometrically. Therefore the $\Gamma_{k-1}(\rchi_{k(+)}-\rchi_{k(-)})$ terms of the summation in the forward pass update equation \eqref{smoothed_state_calculation_simplified} will also tend to $\vec{0}$ and the difference in marginal control contribution $\|u_0^{[N]} - u_0^{[N-1]}\|$ from additional prediction horizon length vanishes.
\end{proof}

\begin{corollary} \label{cor: delta p0 -> 0}
    Let Assumption \ref{assumption_kf_converges} hold. The difference in the initial-time smoothed covariances between adjacent prediction horizon lengths $\Psi$, given by
    \begin{equation}
        \label{psi_delta_p_0_s}
        \Psi_N = P_{0[s]}^{[N]} - P_{0[s]}^{[N-1]}
    \end{equation}
    converges to $\mathbb{O}$ geometrically as $N\rightarrow \infty$.
\end{corollary}

\begin{proof}
    In a similar argument to that of Corollary \ref{cor: delta epsilon}, the convergence of the gain term $\Gamma_{N-1} \rightarrow \mathbb{O}$ ensures that the individual terms $\Gamma_{k-1}(P_{k(-)} - P_{k(+)})\Gamma_{k-1}^\top$ in the summation of the smoother covariance equation \eqref{smoother_cov_forward_simplified} will geometrically converge to the infinite-horizon value of $P_{0[s]}^{[\infty]}$. Therefore, $\Psi_N \rightarrow \mathbb{O}$ geometrically as $N \rightarrow \infty$.
\end{proof}

\begin{remark}
    If the system pair $(A, B)$ is uncontrollable, $P_{0[s]}$ will not converge in the undetectable states as $N \rightarrow \infty$, due to the presence of eigenvalues on the unit circle. However, the difference $\Psi$ will still always converge to $\mathbb{O}$.
\end{remark}

\subsection{Early Termination Criteria} \label{subsec:oc_early_term_crit}

Algorithm \ref{alg:oc_foward_only} computes the control action $u_0$ using a forward-in-time recursion and Corollaries \ref{cor: delta epsilon}--\ref{cor: delta p0 -> 0} demonstrate when future information can no longer affect $u_0$ and $P_{0[s]}$. As such, early termination metrics can determine when the prediction horizon has reached a sufficient length. The first such metric is provided by Corollary \ref{cor: delta epsilon} which establishes a $\rho$ as
\begin{equation}
    \label{foward_only_rho_term}
    \rho_k =\| \Gamma_{k-1} K_k \|
\end{equation}
at which point a threshold $\rho \le \delta_1$ sets a minimum look ahead. In practice, the threshold should be chosen based on desired control smoothness, i.e. the maximum $\Delta u_0$ permissible given an unanticipated change in reference at time $N+1$.

Similarly, Corollary \ref{cor: delta p0 -> 0} develops an early termination metric based on the convergence of $P_{0[s]}^{[N]}$, which is analogous to a dual's cost-to-go Riccati equation. Therefore, \eqref{psi_delta_p_0_s} is a valuable metric for determining when $u_0$ has reached an effective infinite-time equivalence and is a first-order optimality condition on $P_{0[s]}$. This informs the $\tau$ termination metric
\begin{align}
    \notag
    \tau_k &= - \text{tr}(\Psi) / \text{~tr} (P_{0[s]}^{[k]})\\
    \label{eqn:explicit_tau_compute_gammas}
    &= - \text{tr}(\Gamma_{k-1} (P_{k(+)} - P_{k(-)}) \Gamma_{k-1}^\top ) / \text{~tr} (P_{0[s]}^{[k]})
\end{align}
where $\tau_k \le \delta_2$ ensures that the algorithm has sufficiently approximated infinite-time regulating control.

Both termination metrics should be enforced in practice because they certify complementary notions of convergence. The metric $\rho_k$ measures the sensitivity of the current control $u_0$ to an additional future residual, and therefore determines whether the available reference sequence has been sufficiently accounted for. In contrast, $\tau_k$ is independent of the reference sequence and measures the remaining marginal reduction in the initial-time smoother covariance, thus verifying convergence towards the infinite-horizon regulating solution.

\subsection{Efficient Single-Inverse Observed Control} \label{subsec:oc_numerical_improvements}

While Algorithms \ref{alg:oc_plain}--\ref{alg:oc_foward_only} provide for a straightforward implementation, further efficiency and numerical stability can be gained by analyzing the structure of the smoother gain accumulator $\Gamma_k$. As given in \eqref{gain_accumulator_eqn}, the accumulator is intuitively the product of the individual RTS smoother gains $(L_0, \cdots, L_k)$; however, the inverse priors in each $L_i$, $i\in[0,k]$ can be to evaluate and accumulate numerical sensitivity. Expanding the first few terms of the smoother gain accumulator \eqref{gain_accumulator_eqn} gives
\begin{equation}
    \label{expanded_gamma}
    \Gamma_k = \Gamma_{-1} P_{0(+)} \Phi_0^\top P_{1(-)}^{-1}  P_{1(+)} \Phi_1^\top P_{2(-)}^{-1}  \prod_{i=2}^{k} L_i.
\end{equation}
Replacing each $P_{k(+)}$ with the transpose of \eqref{kf:posterior_cov_update} yields
\begin{align}
    \Gamma_k &= \Gamma_{-1} P_{0(-)}\left(\mathbb{I} - K_0 H_0\right)^\top \Phi_0^\top P_{1(-)}^{-1} \notag\\
          & \qquad P_{1(-)}^\top (\mathbb{I} - K_1 H_1)^\top \Phi_1^\top P_{2(-)}^{-1} \prod_{i=2}^{k}L_i \notag\\
    \notag
    &= \Gamma_{-1} P_{0(-)}(\mathbb{I} - K_0 H_0)^\top \\&\qquad\left( \prod_{i=1}^k \Phi_{i-1}^\top (\mathbb{I} - K_i H_i)^\top \right) \Phi_{k}^\top P_{k+1(-)}^{-1}.
\end{align}
Defining $\Phi_{-1} \equiv \mathbb{I}$, the smoother gain accumulator becomes
\begin{align}
    \label{numeric_stable_gain_prod}
    \Gamma_k = \Gamma_{-1} P_{0(-)} \Big( \prod_{i=0}^k \Phi_{i-1}^\top (\mathbb{I} - K_i H_i)^\top \Big) \Phi_{k}^\top P_{k+1(-)}^{-1}.
\end{align}
Expanding the Kalman gain in \eqref{eqn_posterior_prior_difference} with \eqref{kf:gain} yields
\begin{equation}
    \label{eqn_expanded_posterior_prior_diff}
    \rchi_{k(+)} - \rchi_{k(-)} = P_{k(-)} H_k^\top (H_k P_{k(-)} H_k^\top + \grave{R})^{-1} r_k.
\end{equation}
Combining \eqref{eqn_expanded_posterior_prior_diff} with the improved gain accumulator \eqref{numeric_stable_gain_prod} and substituting into equation \eqref{smoothed_state_calculation_simplified} gives the result
\begin{equation}
    \label{numeric_stable_state_acc}
    u_{0[s]} = u_{0(+)} +  \sum_{k=1}^{N} \left[G_{k-1} \Phi_{k-1}^\top S_k r_k\right]
\end{equation}
where $G_k$ and $S_k$ are the inverse free gain and the observed innovation matrices, respectively, given by
\begin{subequations}
    \begin{align}
        \label{eq:Phi minus one definition}
        \Phi_{-1} &\equiv \mathbb{I} \\
        \label{eq:G minus one definition}
        G_{-1} &\equiv \begin{bmatrix} \mathbb{O} & \mathbb{I} \end{bmatrix} P_{0(-)} \\ 
        \label{eqn_G_k_def}
        G_k &= \begin{bmatrix} \mathbb{O} & \mathbb{I} \end{bmatrix} P_{0(-)} \prod_{i=0}^k \Phi_{i-1}^\top (\mathbb{I} - K_i H_i)^\top \\
        \label{num_stable_gain_acc}
        G_k &= G_{k-1}\Phi_{k-1}^\top\left(\mathbb{I}-K_k H_k\right)^\top \\
        \label{S_k_definition}
        S_k &= H_k^\top (H_k P_{k(-)} H_k^\top + \grave{R})^{-1}.
    \end{align}
\end{subequations}
Although this reformulation appears messier, it requires only one inverse per time horizon update in $S_k$, which is guaranteed to be well-defined by the positive-definite $\grave{R}$. This single inverse can be reused in the Kalman gain calculation, giving 
\begin{equation}
    \label{eqn_new_S_based_K}
    K_k=P_{k(-)} S_k.
\end{equation}
As in \eqref{smoothed_state_calculation_simplified}, all required information can be accumulated from the forward filter pass, minimizing computational burden.

\begin{algorithm}[!t]
    \caption{Efficient Observed Control}\label{alg:efficient_update}
    \begin{algorithmic}[1]
    \State {\textbf{initialize} $\rchi_{0(-)} \gets \begin{bmatrix}
        \hat{x} & u_\text{last}
    \end{bmatrix}^\top $, $P_{0(-)} \gets \grave{Q}$,} 
    \State{$\text{tr}(P_{0[s]}) \gets \text{tr}(P_{0(-)})$, $u_{0[s]} \gets u_\text{last}$, $G_{-1}\gets \begin{bmatrix} \mathbb{O} & \mathbb{I} \end{bmatrix}\grave{Q}$}
    \\
    \For {$k \in [0, \cdots, N]$} \label{alg3_line:for_loop_start} \Comment{Single Forward Pass}
        \State {\textbf{perform} Algorithm \ref{alg:oc_plain} \textbf{lines} \ref{alg1_line:kf_start}-\ref{alg1_line:kf_residual}}
        \State {$S_k \gets H_k^\top [H_k P_{k(-)} H_k^\top + \grave{R}_k]^{-1}$} \Comment{One Inverse \eqref{S_k_definition}}
        \State {$K_k \gets P_{k(-)} S_k$} \Comment{\eqref{eqn_new_S_based_K}}
        \State {\textbf{perform} Algorithm \ref{alg:oc_plain} \textbf{lines} \ref{alg1_line:kf_posterior}-\ref{alg1_line:kf_end}}
        \\
        \State {$\Psi_k \gets G_{k-1} \Phi_{k-1}^\top S_k H_k \Phi_{k-1} G_{k-1}^\top $} \Comment{$\Delta P_{0[s]}$ }
        \State {$u_{0[s]} \gets u_{0[s]} + G_{k-1} \Phi_{k-1}^\top S_k r_k$} \Comment{\eqref{numeric_stable_state_acc}}
        \State{$\text{tr}(P_{0[s]}) \leftarrow  \text{tr}(P_{0[s]}) - \text{tr}\left( \Psi_k \right)$} \Comment{\eqref{numeric_stable_cov_acc}}
        \\
        \State{$\rho_k \gets \|  G_{k-1} \Phi_{k-1}^\top S_k \|$} \Comment{\eqref{eqn_efficient_rho_calc}}
        \State{$\tau_k \gets \text{tr}( G_{k-1} \Phi_{k-1}^\top S_k H_k \Phi_{k-1}G_{k-1}^\top ) / \text{tr} (P_{0[s]}) $} \Comment{\eqref{eqn_efficient_tau_calc}}
        \State{\textbf{if} $\rho \le \delta_1$ \textbf{and} $\tau_k \le \delta_2$ \textbf{break}} \Comment{Early Termination}
        \\
        \State {$G_k=G_{k-1} \Phi_{k-1}^\top (\mathbb{I}-K_k H_k)^\top $} \Comment{\eqref{num_stable_gain_acc}}
    \EndFor\\
    \Return $u_{0[s]}$ \label{alg3_line:return_statement}
    \end{algorithmic}
\end{algorithm}

A similar result for the forward-in-time update of the smoothed covariance $P_{0[s]}$ can be computed. The posterior-prior covariance difference $P_{k(+)}-P_{k(-)}$ in \eqref{smoother_cov_forward_simplified} becomes
\begin{equation}
    P_{k(+)}-P_{k(-)} = -K_k H_k P_{k(-)}
\end{equation}
which can be expressed in terms of $S_k$ as
\begin{equation}
    \label{priorminuspost_covs}
    P_{k(+)}-P_{k(-)} = -P_{k(-)} S_k H_k P_{k(-)}.
\end{equation}
The priors on either side of \eqref{priorminuspost_covs} cancel with the inverses on \eqref{numeric_stable_gain_prod} when plugged into \eqref{smoother_cov_forward_simplified} yielding
\begin{equation}
    \label{numeric_stable_cov_acc}
    P_{0[s]}  = P_{0(+)} - \sum_{k=1}^{N} G_{k-1} \Phi_{k-1}^\top S_k H_k \Phi_{k-1}G_{k-1}^\top
\end{equation}
as the improved forward update equation for the smoothed covariance $P_{0[s]}$. Finally, the $\rho_k$ and $\tau_k$ termination criteria can now be re-considered in terms of $G_{k-1}$ and $S_k$ giving 
\begin{align}
    \label{eqn_efficient_rho_calc}
    \rho_k &= \|  G_{k-1} \Phi_{k-1}^\top S_k \| \\
    \label{eqn_efficient_tau_calc}
    \tau_k &= -\text{tr}( G_{k-1} \Phi_{k-1}^\top S_k H_k \Phi_{k-1}G_{k-1}^\top ) / \text{~tr} (P_{0[s]}^{[k]}) .
\end{align}

This development directly leads to Algorithm \ref{alg:efficient_update}---\textit{Efficient Observed Control}. The algorithm has a very similar computational burden to that of an EKF, with only one matrix inverse (guaranteed positive-definite symmetric) required each step. As in Algorithm \ref{alg:oc_foward_only}, its forward-only structure also affords early termination when a sufficient horizon length is reached.

\subsection{Separability into Reactive and Anticipatory} \label{subsec:oc_anytime}

For linear-quadratic formulations such as \eqref{eqn:lqr_forulation}, Algorithm \ref{alg:efficient_update} shows the recursive nature of the MPC update, which is linear in both the reference trajectory and the current state.

\begin{assumption}
    \label{assumption_full_state_feedback}
    The controller, as implemented in a control loop, has full-state information at the control update time.
\end{assumption}

\begin{theorem}
    \label{thm: klqr and feed forward}
    Let Assumptions \ref{assumpion_reasonable_weights_detectable_sys}--\ref{assumption_full_state_feedback} hold. The optimal predictive control policy can be factored into purely reactive and feed-forward anticipatory components
    \begin{align}
        \label{factored_anytime_anyhorizon_formulation}
        u_0 = u_{\text{last}} + K_{\text{lqr}}r_{0} + \Lambda_{0}z_{0} + \{\cdots\} + \Lambda_{N}z_{N}
    \end{align}
    where the gain matrices are given by:
    \begin{enumerate}
        \item {
            The sequence 
            \begin{align}
                \label{eq: klqr sequence}
                K_{\text{eff}} = G_{-1} S_{0} + \sum_{k=1}^{N} G_{k-1} \Phi_{k-1}^{\top} S_{k} H_{k} \prod_{i=0}^{k-1} \Phi_{i} \left( \mathbb{I} - K_{i}H_{i} \right)
            \end{align}
            converges to the stabilizing LQR gain $K_{\text{lqr}}$ as $N \rightarrow \infty$.
        }
        \item {
            The sequence for $\Lambda_0$ given by
            \begin{align}
                \label{eq: lambda_0}
                \Lambda_{0} &= -\sum_{k=1}^{N}G_{k-1}\Phi_{k-1}^{\top}S_{k}H_{k}\Biggl(\prod_{i=1}^{k-1}\Phi_{i}\left(\mathbb{I}-K_{i}H_{i}\right)\Biggr)\Phi_{0}
            \end{align}
        }
        \item {
            and the sequences for $\Lambda_1, \cdots, \Lambda_N$ given by
            \begin{alignat}{2}
                \label{eq: lambda_k}
                \Lambda_{k} &= G_{k-1}\Phi_{k-1}^{\top}S_{k}  &&-\sum_{j=k+1}^{N} G_{j-1} \Phi_{j-1}^{\top}S_{j}H_{j} \\
                \notag & && \times\Biggl(\prod_{i=k+1}^{j-1}\Phi_{i}\left(\mathbb{I}-K_{i}H_{i}\right)\Biggr)\Phi_{k}K_{k}
            \end{alignat}
            converge to the optimal feed-forward anticipatory gains $(\Lambda_0, \Lambda_1 \cdots,\Lambda_N)$, such that the control contribution $u_{a,n}$ resulting from the n\textsuperscript{th} target $z_n$ is given by $u_{a,n} = \Lambda_{n}z_n$.
        }
    \end{enumerate}
\end{theorem}
\begin{proof}
    The forward dynamics of the prior states for the observed system can be constructed from \eqref{eq: forward pass} as
    \begin{subequations}
        \begin{alignat}{2}
            x_{0(-)} & = && -r_0 + z_0\\
            x_{1(-)} & = && -\Phi_{0}\left(\mathbb{I}-K_{0}H_{0}\right)r_{0}+\Phi_{0}z_{0}\\
            x_{2(-)} & = && -\Phi_{1}\left(\mathbb{I}-K_{1}H_{1}\right)\Phi_{0}\left(\mathbb{I}-K_{0}H_{0}\right)r_{0}\\
            \notag
                & && + \Phi_{1}\left(\mathbb{I}-K_{1}H_{1}\right)\Phi_{0}z_{0}+\Phi_{1}K_{1}z_{1}\\
            \label{efficient_sequence_of_priors}
            x_{k(-)} & = && -\Biggl(\prod_{i=0}^{k-1}\Phi_{i}\left(\mathbb{I}-K_{i}H_{i}\right)\Biggr)r_{0} \\
            \notag
                & && +\Biggl(\prod_{i=1}^{k-1}\Phi_{i}\left(\mathbb{I}-K_{i}H_{i}\right)\Biggr)\Phi_{0}z_{0}\\
            \notag
                & && +\sum_{j=1}^{k-1}\Biggl(\prod_{i=j+1}^{k-1}\Phi_{i}\left(\mathbb{I}-K_{i}H_{i}\right)\Biggr)\Phi_{j}K_{j}z_{j}.
        \end{alignat}
    \end{subequations}
    This then yields the sequence of residuals for $k = \{1, \cdots, N\}$
    \begin{alignat}{2}
        \label{efficient_sequence_of_residuals}
        r_{k} & = z_{k} && + H_k \left(\prod_{i=0}^{k-1}\Phi_{i}\left(\mathbb{I}-K_{i}H_{i}\right)\right)r_{0}\\
        \notag
            & && - H_{k}\left(\prod_{i=1}^{k-1}\Phi_{i}\left(\mathbb{I}-K_{i}H_{i}\right)\right)\Phi_{0}z_{0}\\
        \notag
            & && -H_{k}\sum_{j=1}^{k-1}\left(\prod_{i=j+1}^{k-1}\Phi_{i}\left(\mathbb{I}-K_{i}H_{i}\right)\right)\Phi_{j}K_{j}z_{j}.
    \end{alignat}
    Combining with \eqref{numeric_stable_state_acc}, the expanded control is given by
    \begin{alignat}{2}
        \label{eqn_expanded_anytime_control_summation}
        u_{0[s]} & = u_{0(+)} && + \sum_{k=1}^{N}G_{k-1}\Phi_{k-1}^{\top}S_{k} \Bigg\{ z_{k} \\
        \notag
            & && + H_{k}\Biggl(\prod_{i=0}^{k-1}\Phi_{i}\left(\mathbb{I}-K_{i}H_{i}\right)\Biggr)r_{0} \\
        \notag
            & && - H_{k}\Biggl(\prod_{i=1}^{k-1}\Phi_{i}\left(\mathbb{I}-K_{i}H_{i}\right)\Biggr)\Phi_{0}z_{0} \\
        \notag 
            & && - H_{k}\sum_{j=1}^{k-1}\Biggl(\prod_{i=j+1}^{k-1}\Phi_{i}\left(\mathbb{I}-K_{i}H_{i}\right)\Biggr)\Phi_{j}K_{j}z_{j} \Bigg\}.
    \end{alignat}
    Using the relation $u_{0(+)} = u_{0(-)} + G_{-1} S_0 r_0$ and rearranging the last summation of \eqref{eqn_expanded_anytime_control_summation} to group the the $z_k$ terms yields
    \begin{alignat}{2}
        \label{eq: anytime sequence expanded}
        u_{0[s]} & = && ~ u_{\text{last}} + G_{-1}S_{0}r_0 \\
        \notag
            & && + \sum_{k=1}^{N} G_{k-1} \Phi_{k-1}^{\top} S_{k} H_{k} \Biggl( \prod_{i=0}^{k-1} \Phi_{i} \left( \mathbb{I} - K_{i}H_{i} \right) \Biggr) r_0 \\
        \notag
            & && - \sum_{k=1}^{N}G_{k-1}\Phi_{k-1}^{\top}S_{k}H_{k}\Biggl(\prod_{i=1}^{k-1}\Phi_{i}\left( \mathbb{I} - K_{i}H_{i}\right)\Biggr)\Phi_{0} z_0 \\
        \notag
            & && + \sum_{k=1}^{N} \Biggl[ G_{k-1}\Phi_{k-1}^{\top}S_{k} - \sum_{j=k+1}^{N} G_{j-1} \Phi_{j-1}^{\top}S_{j}H_{j}\\
        \notag
            & && \qquad \qquad \quad \times \Biggl(\prod_{i=k+1}^{j-1}\Phi_{i}\left(\mathbb{I}-K_{i}H_{i}\right)\Biggr)\Phi_{k}K_{k} \Biggr] z_k
    \end{alignat}
    where correspondence to the theorem subparts is directly observed. For regulators, all trajectory targets $(z_0, \cdots, z_N)$ are trivially zero, leaving only the reactive components multiplied by $r_0$ in \eqref{eq: anytime sequence expanded}. By Theorem \ref{thm:oc_lqr_equilvalence}, this zero-reference arrangement provides the optimal solution to the LQR regulation problem, which for the augmented system has the solution
    \begin{equation}
        u = u_{\text{last}} + K_{\text{lqr}}r_0
    \end{equation}
    where $u_{\text{last}} = u_{0(-)}$. Therefore, the $K_{\text{eff}}$ sequence multiplying $r_0$ in \eqref{eq: anytime sequence expanded} must converge to the infinite horizon LQR gain as $N\rightarrow\infty$. Further, given the equivalence between this control solution and the LQR problem in Theorem \ref{thm:oc_lqr_equilvalence}, the $(\Lambda_0, \cdots, \Lambda_N)$ sequences must also converge to their optimal feed-forward anticipatory trajectory gains.
\end{proof}

\begin{corollary}
    \label{cor: delta_K converges monitonically}
    Let Assumptions \ref{assumpion_reasonable_weights_detectable_sys}--\ref{assumption_full_state_feedback} hold. The convergence of \eqref{eq: klqr sequence} to $K_{\text{lqr}}$ is monotonic and geometric in horizon length $N$.
\end{corollary}
\begin{proof}
    For each additional horizon step $N$, the change in $K_{\text{eff}}$ is the difference $ K_{\text{eff}}^{[N]} - K_{\text{eff}}^{[N-1]}$ which is given by
    \begin{equation}
        \Delta K_{\text{eff}} = G_{N-1}\Phi_{N-1}^{\top}S_{N}H_{N}\prod_{i=0}^{N-1}\Phi_{i}\left(\mathbb{I}-K_{i}H_{i}\right)
    \end{equation}
    which after using the definition \eqref{eq:Phi minus one definition} becomes
    \begin{align}
        \label{phi_negative_1_trick}
        \Delta K_{\text{eff}} = G_{N-1}\Phi_{N-1}^{\top}S_{N}H_{N} \Phi_{N-1} \prod_{i=0}^{N-1} \left(\mathbb{I}-K_{i}H_{i}\right) \Phi_{i-1}.
    \end{align}
    Comparing to a similar $P_{0[s]}^{[N]}-P_{0[s]}^{[N-1]}$ treatment of \eqref{numeric_stable_cov_acc} and noting that the product in \eqref{phi_negative_1_trick} is nearly $G_{N-1}^\top$ in \eqref{eqn_G_k_def} yields
    \begin{equation}
        \Delta P_{0[s]} = - \Delta K_{\text{eff}}P_{0(-)}\begin{bmatrix} \mathbb{O} & \mathbb{I} \end{bmatrix}^\top.
    \end{equation}
    Since additional future information can only reduce initial uncertainty, $P_{0[s]}$ must converge monotonically, that is $\Delta P_{0,[s]}$ must be negative semi-definite. Therefore, as $P_{0(-)}$ is positive definite it cannot change the sign of $\Delta K_{\text{eff}}$ and the difference $\|K_{lqr}-K_{\text{eff}}\|$ must converge to zero monotonically.
\end{proof}

\begin{corollary}
    Let Assumptions \ref{assumpion_reasonable_weights_detectable_sys}--\ref{assumption_full_state_feedback} hold. For some constant $\gamma>0$, the error $\|K_{\text{lqr}} - K_{\text{eff}}\| < \gamma \, \text{tr}\left(\Delta P_{0[s]}\right)$ as $N\rightarrow \infty$.
\end{corollary}
\begin{proof}
    In the proof of Corollary \ref{cor: delta_K converges monitonically}, it is shown that $\Delta P_{0[s]} = - \Delta K_{\text{eff}}P_{0(-)}\begin{bmatrix} \mathbb{O} & \mathbb{I} \end{bmatrix}^\top$. Thus, the convergence rates are the same between $K_{\text{eff}}$ and $P_{0[s]}$ within a factor $\gamma$. The rate of convergence is geometric in the square eigenvalues of the closed-loop system $\left(\mathbb{I} - K_i H_i \right) \Phi_{i-1}$, which for large $N$ becomes the closed-loop infinite-time Kalman system $\left(\mathbb{I}-K_\infty H\right)\Phi$, which by Theorem \ref{thm: linf_stable} has all eigenvalues within the unit circle. Thus, the convergence errors are bounded by
    \begin{alignat}{2}
        \label{eq: K_diff sum}
        \|K_{lqr}-K_{\text{eff}}^{[N]}\|  &<  \gamma_k \sum_{n=N}^\infty \lambda_{\text{max}}^{2n}
        &&<\gamma_k \frac{\lambda_{\text{max}}^{2N}}{1-\lambda_{\text{max}}^2} \\
        \label{eq: LQR bounds} 
        \|P_{0[s]}^{[\infty]} - \sum_{n=N}^\infty \Delta P_{0[s],n}\|  &<  \gamma_p \sum_{n=N}^\infty \lambda_{\text{max}}^{2n}
        &&<\gamma_p \frac{\lambda_{\text{max}}^{2N}}{1-\lambda_{\text{max}}^2}.
    \end{alignat}
    where $\lambda_\text{max}$ is the maximum magnitude (controllable) eigenvalue of the closed-loop Kalman system. 
    The error $\|K_{lqr}-K_{\text{eff}}^{[N]}\|$ and initial covariance $P_{0,[s]}$ both converge quadratically with $N$ and differ by the factor $\gamma = \nicefrac{\gamma_k}{\gamma_p}$. Therefore, the early termination metric $\tau =\text{tr}\left({\Delta P_{0[s]}}/{P_{0[s]}}\right)$ is equivalent to $\|K_{lqr}-K_{\text{eff}}^{[N]}\|<\epsilon\|K_{lqr}\|$ for some $\epsilon$.
\end{proof}

\begin{remark}
    Given \eqref{eq: LQR bounds}, an approximation for the horizon length $N$, given a desired convergence tolerance $\epsilon$ such that 
    \begin{equation}
        \epsilon \approx \frac{\|K_{\text{lqr}} - K_\text{eff}\|}{\|K_{\text{lqr}}\|}
    \end{equation}
    is given from the geometric series as
    \begin{equation}
        \label{eqn_required_N_estimate}
        N\approx\frac{\log\left(\epsilon\right)}{2\log\left(\lambda_\text{max}\right)}.
    \end{equation}
    However, since the convergence of $K_\text{eff}$ is not purely geometric until large $N$, this only acts in practice as a heuristic for the termination criteria used in the algorithms and $K_\text{eff}$ convergence to a ratio of $\epsilon$. 
\end{remark}

\begin{corollary}
    \label{cor:smoothness_of_trajectory}
    Let Assumptions \ref{assumpion_reasonable_weights_detectable_sys}--\ref{assumption_full_state_feedback} hold. There exists a finite horizon length $N$ such that the anticipatory effect of including additional information from an infinite horizon Lipschitz continuous reference falls below any desired value $\epsilon$, that is
    \begin{equation}
        \|u^{[\infty]} -u^{[N]}\| < \epsilon \mathrm{~if~} \|z_k -z_{k-1}\| < \delta \, \forall \, k>N.
    \end{equation}
\end{corollary}
\begin{proof}
Consider the any-time any-horizon control factorization \eqref{factored_anytime_anyhorizon_formulation} for an infinite horizon, expressed as
\begin{align}
    u^{[\infty]} &= u_{\text{last}} + K_{\text{lqr}}r_{0} + \Lambda_{0}z_{0}+\Lambda_{1}\left(z_{0}+z_{1}-z_{0}\right)+ \{\cdots\} \nonumber\\ 
    &= u_{\text{last}}+K_{\text{lqr}}r_{0}+\Lambda_{0}z_{0}+\Lambda_{1}\left(z_{0}+\Delta z_{1}\right)+\{\cdots\}\\
    &= u_{\text{last}}+K_{\text{lqr}}r_{0}+\sum_{n=0}^{\infty}\Lambda_{n}z_{0}+\sum_{n=1}^{\infty}\left(\sum_{m=m}^{\infty}\Lambda_{m}\right)\Delta z_{n}. \nonumber
\end{align}
Thus, control difference $u^{[\infty]} - u^{[N]}$ is given by
\begin{equation}
    \label{eqn_proof_infinite_difference}
    u^{[\infty]}-u^{[N]} = \sum_{n=N+1}^{\infty}\left(\sum_{m=m}^{\infty}\Lambda_{m}\right)\Delta z_{n}.
\end{equation}
Applying the triangle inequality to \eqref{eqn_proof_infinite_difference} yields the inequality
\begin{align}
    \label{eqn_proof_normed_inequality}
    \|u^{[\infty]}-u^{[N]}\| & \leq \sum_{n=N+1}^{\infty}\left(\sum_{m=n}^{\infty}\left\Vert \Lambda_{m}\right\Vert \right)\left\Vert \Delta z_{n}\right\Vert .
\end{align}
The norm of each $\Lambda_m$ term, defined in \eqref{eq: lambda_k}, is bounded by
\begin{align}
    \label{eqn_proof_bounded_lamda_norm}
    \| \Lambda_{m} &\|  \le \| G_{m-1} \Phi^\top S_m \|  \\
        & + \sum_{i=m}^{N} \left\Vert G_{i-1} \Phi^{\top} S_{i} H \Phi \left( \prod_{k=m+1}^i \left( \mathbb{I} - K_{k-1} H \right) \Phi \right) K_i \right\Vert. \nonumber
\end{align}
having again applied the triangle inequality. Using the fact that
\begin{align}
    \label{eqn_proof_lambda_max_n}
    \| \Phi_{i-1}^\top (\mathbb{I} - K_iH)^\top \| \le \lambda_{max}
\end{align}
and the multiplicative property of norms, \eqref{eqn_proof_bounded_lamda_norm} becomes
\begin{align}
    \| \Lambda_{n} \|  & \le \gamma_1 \lambda_{max}^{m} + \sum_{i=m}^{\infty} (\lambda_{max}^{i}) \gamma_2 (\lambda_{max}^{i-m})
\end{align}
noting the presence of \eqref{eqn_proof_lambda_max_n} in the $G$ terms given by \eqref{eqn_G_k_def} and defining $\gamma_1 = \| \Phi^\top S_n \|$ and $\gamma_2 = \| \Phi^{\top} S_{i} H \Phi \|\|K_i\|$. Thus,
\begin{align}
    \label{eqn_proof_summation_change}
    \| \Lambda_m \|&\le \gamma_1 \lambda_{max}^{m} + \gamma_2 \lambda_{max}^{2m} \sum_{i=0}^{\infty} \lambda_{max}^{2i} \\
    \label{eqn_proof_result1}
    &\le \gamma_1 \lambda_{max}^{m} + \frac{\gamma_2 \lambda_{max}^{2m}}{1 - \lambda_{max}^{2}}
\end{align}
noting the change in the summation range of \eqref{eqn_proof_summation_change}. Plugging the result \eqref{eqn_proof_result1} into \eqref{eqn_proof_normed_inequality} and applying the triangle inequality and summation property used in \eqref{eqn_proof_summation_change} yields the result
\begin{equation}
    \label{corollary_3_4_result}
    \|u^{[\infty]}-u^{[N]}\|\le \left(
        \frac{\gamma_1 \lambda_{max}^{N+1}}{(1-\lambda_{max})^2}
        + \frac{\gamma_2 \lambda_{max}^{2(N+1)}}{(1-\lambda_{max}^2)^3}
    \right) \delta
\end{equation}
which completes the proof as the terms of \eqref{corollary_3_4_result} will approach $0$ as $N\rightarrow \infty$ since $\lambda_{max}$ is guaranteed to be within the unit circle $(| \lambda_{max} | < 1)$ by Theorem \ref{thm: linf_stable} for a controllable system. However, if the system passes through an uncontrollability, the eigenvalues appear on the unit circle $(|\lambda_{max}| = 1)$ for those horizon steps, and additional look-ahead is required.
\end{proof}

\begin{theorem}
    \label{thm:anytime_separable_anticip}
    Let Assumptions \ref{assumpion_reasonable_weights_detectable_sys}--\ref{assumption_full_state_feedback} hold.
    The anticipatory control contribution of a prediction horizon $u_a$ can be computed independently of the current state based on the current and upcoming references $(z_0, \cdots, z_N)$. This anticipatory component can be added to a purely reactive $K_{\text{lqr}}$ control law, enabling any-time any-horizon MPC and computed using the same Kalman filtering recursions. This yields the control policy
    \begin{equation}
        \label{desired_anytime_policy}
        u_0 = K_{\text{lqr}} r_0 + \sum_{k=0}^{N}u_{a,k}
    \end{equation}
    which enables any-time termination of the algorithm, while ensuring at minimum LQR's stability and performance.
\end{theorem}

\begin{proof}
    To compute the anticipatory control contribution using only $(z_0, \cdots, z_N)$, initialize the Kalman filter with $x_{0(-)}=z_0$. Then the first residual is $r_0 = 0$, which removes the current state-error feedback already handled by $K_{\text{lqr}}$ and yields the sequence of $z_0$-relative residuals given by 
    \begin{alignat}{2}
        \label{eqn_anytime_residuals}
        r_{k} & = z_{k} &&  - H_{k}\left(\prod_{i=1}^{k-1}\Phi_{i}\left(\mathbb{I}-K_{i}H_{i}\right)\right)\Phi_{0}z_{0}\\
        \notag
            & && -H_{k}\sum_{j=1}^{k-1}\left(\prod_{i=j+1}^{k-1}\Phi_{i}\left(\mathbb{I}-K_{i}H_{i}\right)\right)\Phi_{j}K_{j}z_{j}
    \end{alignat}
    for $k\in[1,N]$. Combining \eqref{eqn_anytime_residuals} with the efficient algorithm's summation \eqref{numeric_stable_state_acc} yields the anticipatory control contribution identical to that of \eqref{eq: anytime sequence expanded} in the proof of Theorem \ref{thm: klqr and feed forward} without the terms multiplied by $r_0$. 
    Therefore, with the initialization of $x_{0(-)} = z_0$, the same forward-only approach as Algorithm \ref{alg:efficient_update} computes the anticipatory control contribution sum $u_a$ for any horizon $N$ independent of the current system state error.
\end{proof}

Algorithm \ref{alg:anytime-any-horizon MPC}---\textit{Any-Time Any-Horizon Observed Control}---follows directly from Theorem \ref{thm:anytime_separable_anticip}, which establishes that anticipatory control contributions can be computed independently of the current state. It parallels Algorithm \ref{alg:efficient_update}, but with the key initialization difference of $x_{0(-)}=z_0$. Separately, the infinite-time $K_{\text{lqr}}$ gain is applied to the current error, ensuring at worst, the stability and performance of LQR. Furthermore, the forward pass may be terminated at any point while preserving these guarantees, thereby motivating the algorithm's name.

\begin{remark}
    Theorem \ref{thm: klqr and feed forward} proves that linear-quadratic MPC is the superposition of a reactive and anticipatory component.
    Theorem \ref{thm:anytime_separable_anticip} proves that anticipatory controls $u_a$ do not require any state information and can therefore be precomputed and tabulated for the entire trajectory using Algorithm \ref{alg:anytime-any-horizon MPC}.
\end{remark}

\section{Nonlinear Extensions}
\label{sec:nonlinear_systems_extensions}

Observed control readily extends to MPC formulations with nonlinear dynamics and nonquadratic objectives. Consider the nonlinear discrete-time plant given by
\begin{align}
    x_k &= g(x_{k-1}, u_{k-1})
\end{align}
and define the augmented state $\rchi_k$ as before \eqref{augmented_system}. The nonlinear augmented state transition function $f(\rchi_k)$ is defined as
\begin{align}
    f(\rchi_{k-1}) = \begin{bmatrix}
        g(x_{k-1}, u_{k-1}) \\
        u_{k-1}
    \end{bmatrix}.
\end{align}
The corresponding nonlinear filtering model is 
\begin{subequations}
    \label{nonlinear_augmented_filtering_model}
    \begin{align}
        \label{nonlinear_augmented_filtering_model:dynamics}
        \rchi_{k} &= f(\rchi_{k-1}) + \grave{w} \qquad & \grave{w} \sim \mathcal{N}(0, \grave{Q}) \\
        \label{nonlinear_augmented_filtering_model:measurements}
        z_k &= h (\rchi_k) + \grave{v}            & \grave{v} \sim \mathcal{N}(0, \grave{R})
    \end{align}
\end{subequations}
and analogous to \eqref{eqn:eks_minimization}, the negative log-likelihood of the nonlinear maximum a posteriori, is the smoother objective
\label{eqn:eks_nonlinear_minimization}
\begin{alignat}{2}
    \notag
    J(\rchi_0, \cdots, \rchi_N) & =\: && \|(\rchi_0-\bar{\rchi}_0)\|_{P_{0}^{-1}}^2 + \sum_{k=0}^N \|z_k - h(\rchi_k)\|_{R^{-1}}^2 \\
    & && + \sum_{k=1}^{N} \|\rchi_{k} - f(\rchi_{k-1}) \|_{Q^{-1}}^2.
    \label{eqn:nonlinear_map_obj}
\end{alignat}
Following the same measurement and covariance construction as Theorem \ref{thm:oc_lqr_equilvalence} for  \eqref{eqn:nonlinear_map_obj} produces the exact nonlinear-quadratic dual MPC reference tracking formulation. For general nonlinear objectives, the construction of the appropriate measurement model ($z_k$, $H_k$, and $R_k$) is covered in Section \ref{subsec:oc_nonquadratic_objfun_extensions}.

\begin{algorithm}[!t]
    \caption{Any-Time Any-Horizon Observed Control}\label{alg:anytime-any-horizon MPC}
    \begin{algorithmic}[1]
  
    \State{\textbf{initialize} $P_{0(-)} \gets \grave{Q}$, $G_{-1}\gets \begin{bmatrix} \mathbb{O} & \mathbb{I} \end{bmatrix}\grave{Q}$}
    \State {$\text{tr}(P_{0[s]}) \gets \text{tr}(P_{0(-)})$, $\rchi_{0(-)} \gets z_0$, $u_{a} \gets \mathbf{0}$}
    \\
    \State {$u_0[s] \gets$ \textbf{perform} Algorithm \ref{alg:efficient_update} \textbf{lines} \ref{alg3_line:for_loop_start}-\ref{alg3_line:return_statement}}
    \State{}\Comment{Anticipatory component $u_{a}$ is equivalent to $u_{0[s]}$}\\
    \\
    \Return $u_\text{last} + K_\text{lqr}\left(z_0 - H_0 \rchi_0\right)+u_{0[s]}$
    \end{algorithmic}
\end{algorithm}

\subsection{Iterated Extended RTS Smoothing}

A natural nonlinear backend for Algorithm \ref{alg:oc_plain} is the iterated extended Kalman filter and RTS smoother (IERTS); for further detail, the reader is referred to \cite{book_sarkka_2023_bayesian_filtering_and_smoothing}. Given a nominal trajectory, denoted $(\widetilde{\rchi}_0, \cdots, \widetilde{\rchi}_N)$, equations \eqref{nonlinear_augmented_filtering_model:dynamics}--\eqref{nonlinear_augmented_filtering_model:measurements} are linearized
\begin{subequations}
    \begin{align}
        f(\rchi_{k-1}) & \approx f(\widetilde{\rchi}_{k-1}) + \widetilde{F}_{k-1} (\rchi_{k-1} - \widetilde{\rchi}_{k-1}) \\
        h(\rchi_{k}) & \approx h(\widetilde{\rchi}_{k}) + \widetilde{H}_{k} (\rchi_{k} - \widetilde{\rchi}_{k}).
    \end{align}
\end{subequations}
Where the sensitivities $\widetilde{F}_{k-1}$ and $\widetilde{H}_{k}$ are defined as
\begin{align}
    \label{linearized_sensitivities}
    \widetilde{F}_{k-1} = \frac{\partial f(\widetilde{\rchi}_{k-1})}{\partial \widetilde{\rchi}_{k-1}}
    &&
    \widetilde{H}_{k} = \frac{\partial h(\widetilde{\rchi}_{k})}{\partial \widetilde{\rchi}_{k}}.
\end{align}
The Kalman filtering forward update equations \eqref{kf:dynamics}--\eqref{kf:gain} become their time-indexed counterparts where the sensitivities \eqref{linearized_sensitivities} are used instead of the normal EKF linearizations. The prior prediction \eqref{kf:dynamics} and residual \eqref{kf:anticipated_residual} equations become
\begin{subequations}
    \begin{align}
        \rchi_{k(-)} &= f(\widetilde{\rchi}_{k-1}) + \widetilde{F}_{k-1} (\rchi_{k-1(+)} - \widetilde{\rchi}_{k-1}) \\
        r_k &= z_k - h(\widetilde{\rchi}_k ) - \widetilde{H}_{k} (\rchi_{k(-)} - \widetilde{\rchi}_{k}).
    \end{align}
\end{subequations}
Similarly, the RTS smoother gain in equations \eqref{eq:kalman_smoother_state}--\eqref{eq:kalman_smoother_covariance} becomes a version constructed about the linearized trajectory
\begin{align}
    \widetilde{L}_k &= P_{k(+)} \widetilde{F}_{k}^\top P_{k+1(-)}^{-1}.
\end{align}
The IERTS smoother may be initialized with a single RTS pass and repeating $(\widetilde{\rchi}_0, \cdots, \widetilde{\rchi}_N) \gets (\rchi_{0[s]}, \cdots, \rchi_{N[s]})$ between passes yields the standard algorithm, which is equivalent to applying Gauss-Newton iterations to \eqref{eqn:nonlinear_map_obj} \cite{book_sarkka_2023_bayesian_filtering_and_smoothing, bell_1994_iterated_kalman_smoother_as_gauss_newton_method}.

With IERTS, observed control inherits a Gauss-Newton interpretation. Therefore, IERTS-based observed control is a local nonlinear MPC method, analogous in this respect to SQP, RTI, iLQR, and DDP. A key difference, however, is that observed control's first forward filter pass converges towards the reference trajectory - thus producing high-quality linearizations for the smoother pass. When the smoother iterations converge to a fixed point satisfying the standard second-order sufficient conditions, the resulting smoother states contain the controls that locally minimize the dual MPC problem. Additionally, for challenging problems, Levenberg--Marquardt regularization can improve convergence behavior \cite{sarkka_2020_levenberg_marquardt_ieks}.

More generally, the observed control framework is backend agnostic at the modeling level. Any recursive filter-smoother combination that remains well posed under the deterministic process model $\grave{Q}$ and initial covariance $P_0$ can be used. The resulting controller inherits the same approximation quality, convergence region, and computational performance as the chosen backend has the corresponding nonlinear estimation problem. For example with UKF, observed control inherites a derivative free mode, like sampling-based controllers. Thus, IERTS provides Gauss-Newton local optimality properties, while EKF, UKF, cubature, particle, or other nonlinear smoothing backends provide valid approximate observed-control updates whose optimality claims are only as strong as the corresponding estimator approximation.

For nonlinear IERTS, the $\rho_k$ and $\tau_k$ metrics are evaluated about the current iteration's trajectory and therefore should be interpreted only as local indicators of horizon sufficiency until the entire horizon is converged. For approximating extensions such as EKF and UKF, $\rho_k$ and $\tau_k$ metrics should likewise be interpreted only as approximate, local indicators.

\subsection{Real-Time Iterative Hot Starting}

The IERTS backend enables a real-time iterative implementation of observed control (RTI-OC). The recursive smoothing analogue of RTI schemes \cite{li_1989_multistep_newton_type_control_strategy, ohtsuka_2004_gmres_method_for_nl_rhc, diehl_2001_real_time_optimization_for_large_scale_nonlinear_processes_phd, diehl_2005_rti_scheme_for_nonlinear_opt_feedback_ctrl, zavala_2009_advanced_step_controller, houska_2011_autogenerated_rti_nonlinear_mpc_microsecond, houska_2011_acado_toolkit, manual_ACADO_toolkit, verschueren_2018_acados_sofware_paper}, the RTI-OC procedure is as follows. At each control update, the previous update's smoothed augmented trajectory is shifted to the present time and the final stage extended by a dynamics roll out. The initial state is updated with the current plant state and one IERTS iteration is performed. The control contained in the first smoothed state is then applied to the plant, while the remaining smoothed trajectory hot-starts the subsequent update.

Because one IERTS iteration is equivalent to a Gauss--Newton, the nonlinear OCP is not completely solved to convergence during each update. Instead, the shifted solution tracks the moving local minimizer through successive refinements. With one shift per update, future portions of the trajectory undergo many refinements before reaching the beginning of the prediction horizon. Consequently, well-converged control actions are obtained when the OCP varies sufficiently slowly between updates and the shifted trajectory remains within the local convergence region of the Gauss--Newton iterations.

Notably, unlike conventional RTI implementations that use a fixed prediction horizon selected offline, RTI-OC can adapt the horizon online by extending or terminating the forward recursion. 
During the forward filter pass, the horizon can be extended by appending additional measurements until the local values of $\rho_k$ and $\tau_k$ are satisfied. Alternatively, to shorten the horizon, the same  Section \ref{subsec:early_term_discussion} discussion applies.

\begin{table}[t!]
    \centering
    \resizebox{\linewidth}{!}{%
    \begin{tabular}{@{} l @{\hspace{2pt}} r@{\,=\;}l@{}c @{\hspace{8pt}} r@{\,=\;}l@{}c @{\hspace{8pt}} r@{\,=\;}l@{}c  @{}}
    \toprule
        & \multicolumn{3}{@{}p{2.2cm}@{}}{\textbf{Duality-Based \newline Measurement (a)}}
        & \multicolumn{3}{@{}p{2.2cm}@{}}{\textbf{One-Step Newton \newline Measurement (b)}}
        & \multicolumn{3}{@{}p{2.2cm}@{}}{\textbf{Gradient-Based \newline Measurement (c)}}\\ \midrule
    Residual 
        & $r_k$ & $z_k - \hat{\rchi}_k$
        & \eqref{kf:anticipated_residual}
        & $r_k$ & $-\frac{\partial^2 J}{\partial \rchi_k^2}^{\dagger} \frac{\partial J}{\partial \rchi_k}$
        & \eqref{newton_measurement_update_defs}
        & $r_k$ & $\frac{\partial J}{\partial \rchi_k}$ 
        & \eqref{gradient_measurement_defs} \\
    Sensitivity
        & $H_k$ & $\mathbb{I}_\eta$
        & \eqref{eqn:oc_lqr_synth_measurement}
        & $H_k$ & $\frac{\partial^2 J}{\partial \rchi_k^2}^{\dagger} \frac{\partial^2 J}{\partial \rchi_k^2}$
        & \eqref{newton_measurement_update_defs}
        & $H_k$ & $\frac{\partial^2 J}{\partial \rchi_k^2}$ 
        & \eqref{gradient_measurement_defs} \\
    Covariance 
        & $\grave{R}_k$ & $\begin{bsmallmatrix}
            Q & M \\ M^\top & R
        \end{bsmallmatrix}^{-1}$
        & \eqref{eqn:oc_lqr_qr_mapping}
        & $\grave{R}_k$ & $\frac{\partial^2 J}{\partial \rchi_k^2}^{\dagger}$
        & \eqref{newton_measurement_update_defs}
        & $\grave{R}_k$ & $\alpha \mathbb{I}_\eta$ 
        & \eqref{gradient_measurement_defs}  \\
    \bottomrule
    \end{tabular}
    }%
    \caption{
    Duality-informed measurements (a) for exact quadratic objectives. Newton-like (b) and gradient-based (c) observed control measurement modes for nonquadratic objectives.
    \vspace{-1em}}
    \label{table:comp_agg_meas_summary}
\end{table}

\subsection{Extensions to Nonquadratic Objective Functions} \label{subsec:oc_nonquadratic_objfun_extensions}

For objectives of the form \eqref{eqn:disc_time_opt_ctrl_formulation:obj_fn} with nonquadratic stage costs $J_k(x_k, u_k)$, both Newton-like and gradient-based methods are presented. Algorithms \ref{alg:oc_plain}--\ref{alg:efficient_update} readily incorporate these extensions by replacing the corresponding measurement equations with those for the desired measurement mode as summarized in Table \ref{table:comp_agg_meas_summary}.
The Newton-based measurement follows from the standard approach of using sequential Newton-like steps to address nonlinearities in optimization. For standard MPC objective functions, i.e. explicit scalar objectives, this is the preferred method. For problems where an explicit scalar objective function is challenging to define or not practical to tune, the gradient-based measurement is the preferred approach.

\begin{figure*}[t!]
    \centering
    \begin{subfigure}[b]{0.32\textwidth}
        \centering
        \includegraphics[width=\textwidth]{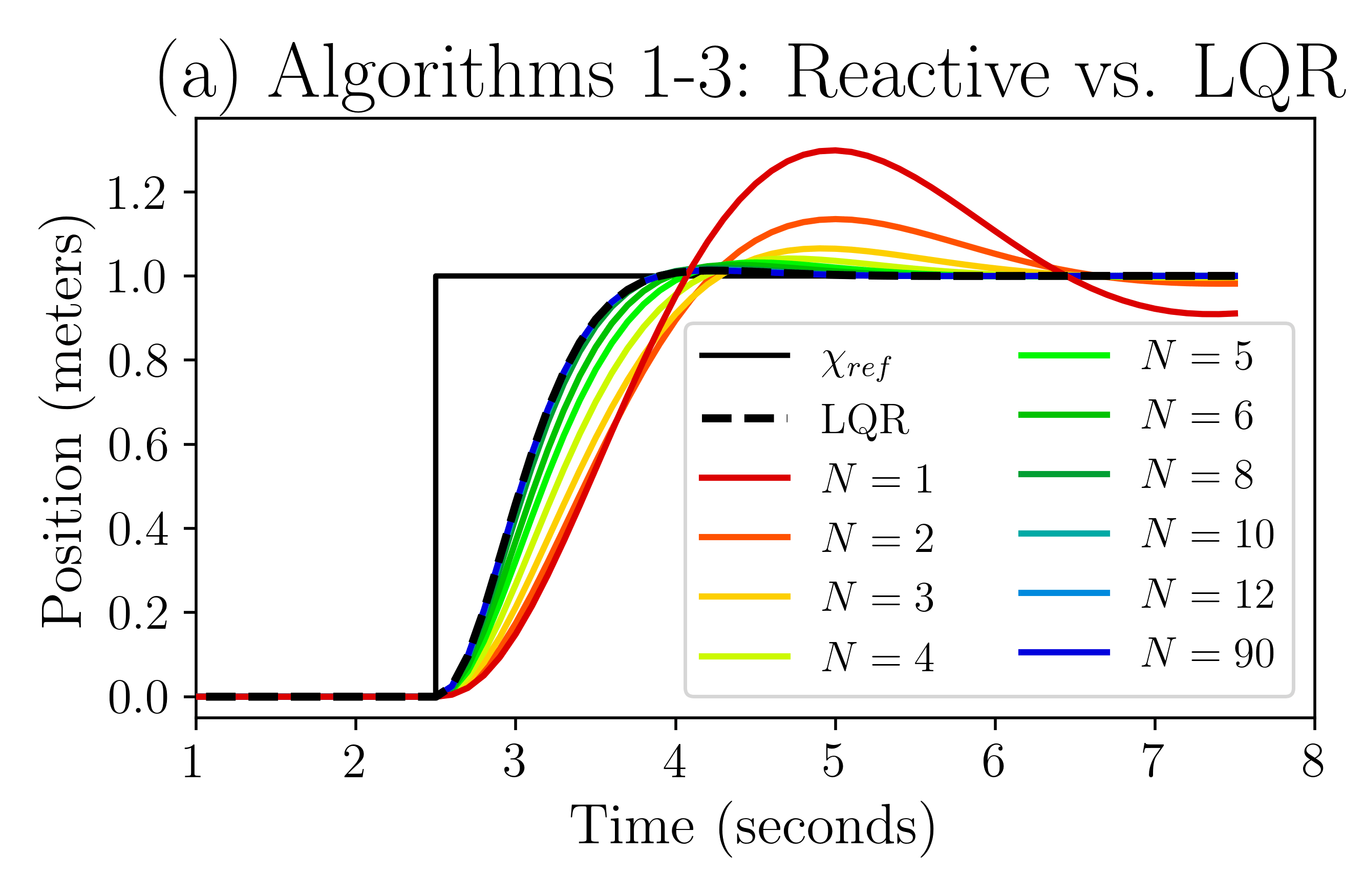}
        \phantomsubcaption%
        \vspace{-1.5em}
        \label{fig:msd_horizon_reactive}
    \end{subfigure}%
    \hfill
    \begin{subfigure}[b]{0.32\textwidth}
        \centering
        \includegraphics[width=\textwidth]{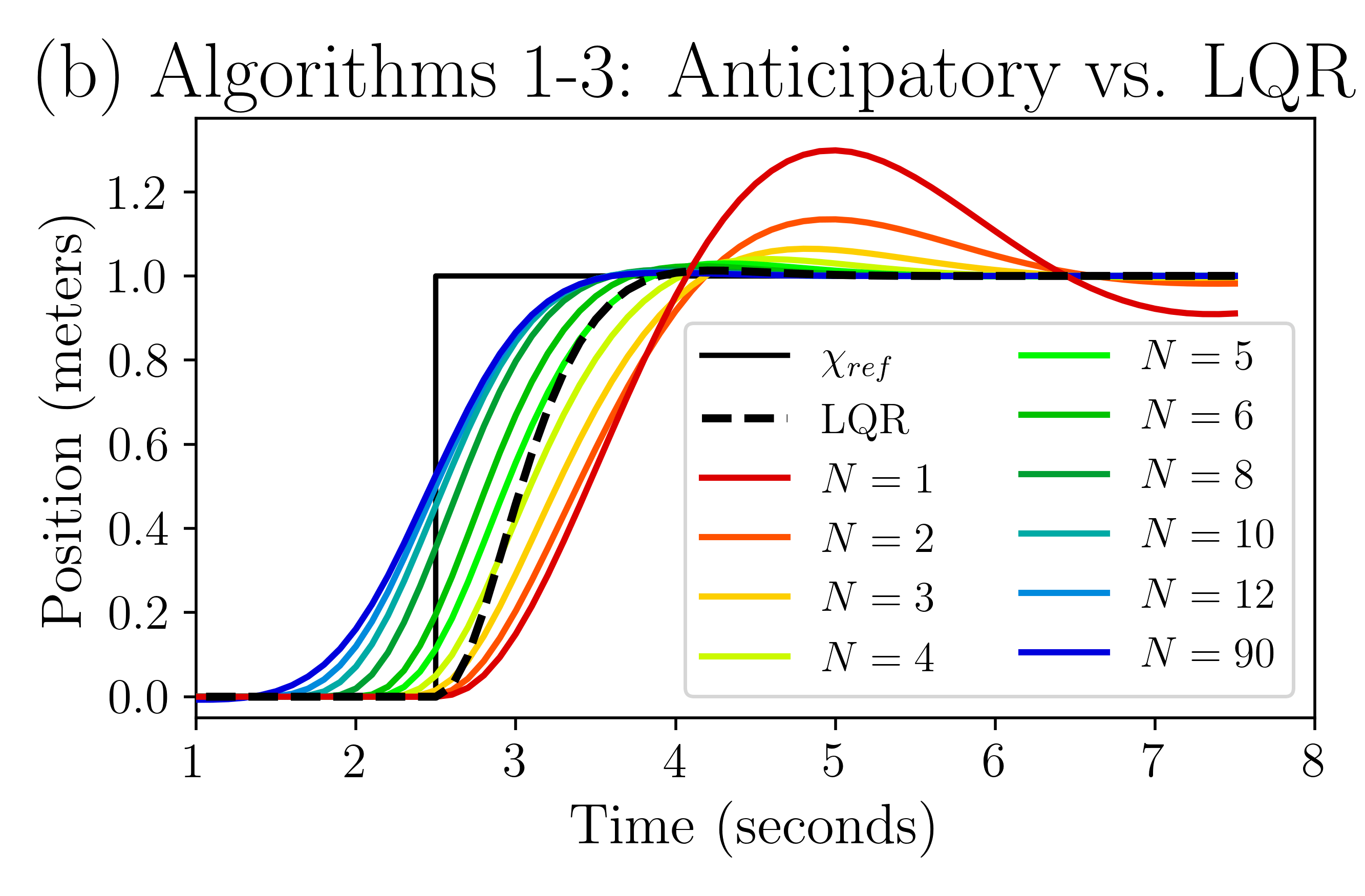}
        \phantomsubcaption%
        \vspace{-1.5em}
        \label{fig:msd_horizon_anticipatory}
    \end{subfigure}%
    \hfill
    \begin{subfigure}[b]{0.32\textwidth}
        \centering
        \includegraphics[width=\textwidth]{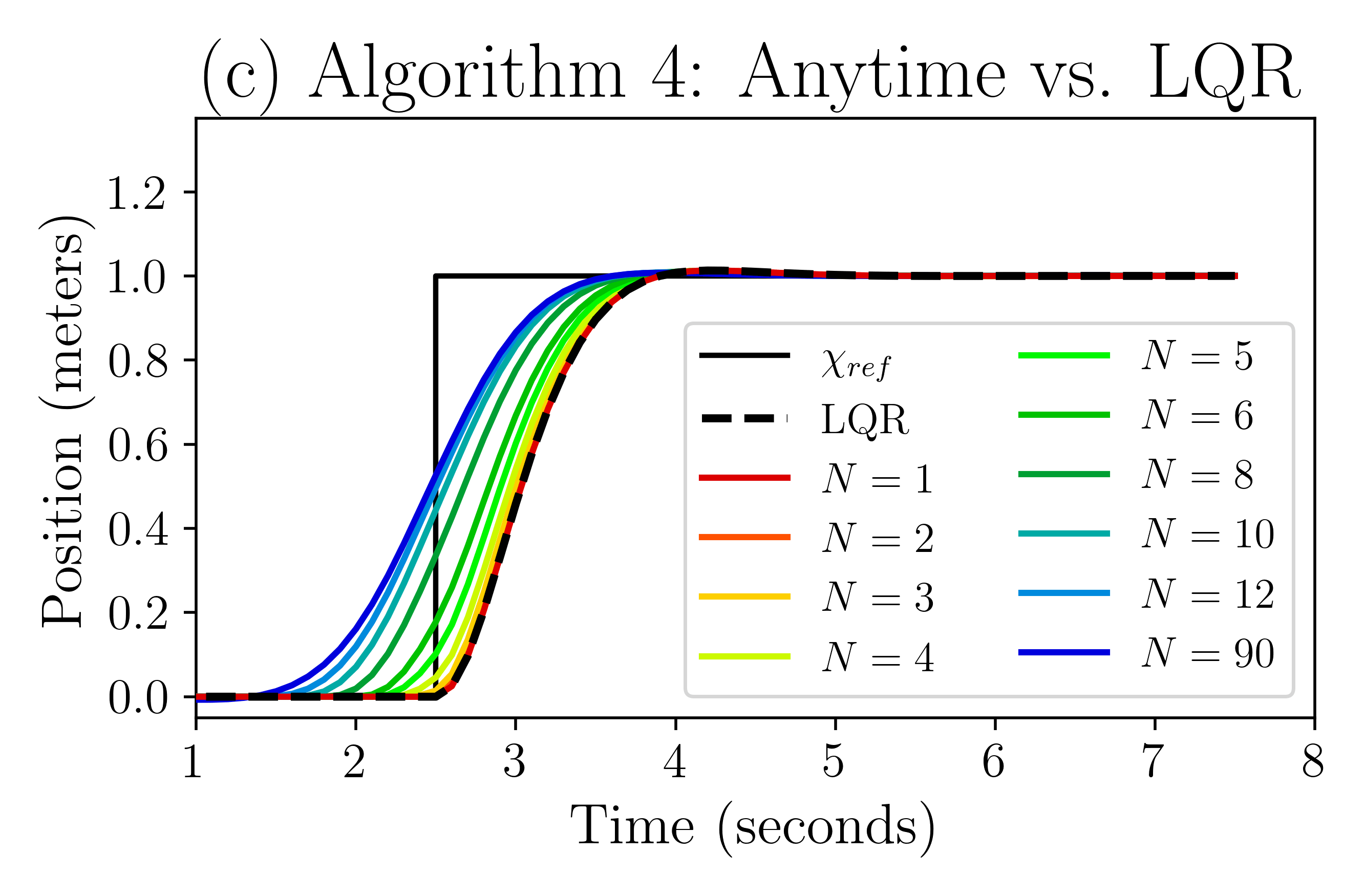}
        \phantomsubcaption%
        \vspace{-1.5em}
        \label{fig:msd_horizon_anytime}
    \end{subfigure}
    \caption{
    Step responses of Algorithms \ref{alg:oc_plain}-\ref{alg:anytime-any-horizon MPC} $(\Delta t = 0.1)$ in reactive (\ref{fig:msd_horizon_reactive}), predictive (\ref{fig:msd_horizon_anticipatory}), and any-time any-horizon (\ref{fig:msd_horizon_anytime}) modes. Without trajectory anticipation, as $N\rightarrow\infty$, Algorithms \ref{alg:oc_plain}-\ref{alg:efficient_update} converge to the optimal infinite-horizon discrete LQR solution; with anticipation, all converge to the predictive optimal; notably, Algorithm \ref{alg:anytime-any-horizon MPC} achieves the reactive optimal with $N=0$.
    \vspace{-1.5em}}
    \label{fig:msd_step_response_trajectories}
\end{figure*}

\subsubsection{Information-Regularized Newton Measurement}

Consider the second order Taylor-series expansion of the single-step stage contribution $J(\rchi_k)$, expressed with the augmented state,
\begin{align}
    \label{eqn_nonlinear_obj_taylor}
    J(\bar{\rchi}_k + \delta\rchi ) \approx J(\bar{\rchi}_k) + J_x ^\top (\delta\rchi) + \frac{1}{2} (\delta\rchi)^\top J_{xx} (\delta\rchi)
\end{align}
where the jacobian and hessian are defined as
\begin{align}
    J_x = \frac{\partial J}{\partial \bar{\rchi}_k}
    &&
    J_{xx} = \frac{\partial^2 J}{\partial \bar{\rchi}_k^2}.
\end{align}
The minimizer to \eqref{eqn_nonlinear_obj_taylor} is the Newton step
\begin{align}
    \label{eqn_newton_step_minimizer}
    \delta\rchi = -\left( J_{xx} \right) ^{-1} J_x.
\end{align}
Theorem \ref{thm:oc_lqr_equilvalence} proves the duality mapping for the pure quadratic objective function; thus, it follows that the quadratic approximation \eqref{eqn_nonlinear_obj_taylor} will follow a similar development in the measurement sensitivity $H_k$ and covariance $\grave{R}$. Taking the newton step \eqref{eqn_newton_step_minimizer} as the residual yields the measurement model
\begin{align}
    \label{newton_measurement_update_defs}
    H_k = J_{xx}^{-1}J_{xx} =\mathbb{I} && \grave{R} = \left( J_{xx} \right)^{-1} && r_k= - J_{xx} ^{-1} J_x.
\end{align}
Expanding the posterior update equation \eqref{kf:posterior_state_update} yields
\begin{align}
    \label{numerically_safe_newton_step}
\rchi_{k(+)} = \rchi_{k(-)} - P_{k(-)}\left(P_{k(-)} + J_{xx}^{-1} \right)^{-1} J_{xx}^{-1} J_{x}
\end{align}
which can be written in the information form as
\begin{align}
    \label{eqn_newton_measurement_information_form}
    \rchi_{k(+)} = \rchi_{k(-)} - \left( P_{k(-)}^{-1} + J_{xx} \right)^{-1} J_x
\end{align}
which reveals the information-regularized Newton-step structure of the posterior Kalman update. In practice, equation \eqref{numerically_safe_newton_step} should be used to avoid inverting singular priors.
Just as in nonlinear SQP solvers, if the Hessian is rank deficient, a symmetry-preserving pseudo-inverse should be used in \eqref{newton_measurement_update_defs}.

\subsubsection{Gradient-Based Measurement} Often it is more natural to express the desired anticipated residual, i.e., regret, $r_k$ directly or heuristically, instead of as an explicit derivative of a scalar function. This is advantageous as it is often easier to define a relative \textit{which way is better} $(\nicefrac{\partial J}{\partial \rchi})$ rather than an \textit{exactly how bad is it} $(J)$. Thus, in general, the residual function $r_h(\rchi_{k(-)})$ can be incorporated in the same guise as the EKF, where
\begin{align}
    \label{gradient_measurement_defs}
    r_k = r_h({\rchi_{k{(-)}}}) &&
    H_k = \frac{\partial r_h({\rchi_{k{(-)}}})}{\partial \rchi_{k(-)}} &&
    \grave{R_k} = \alpha \mathbb{I}.
\end{align}
The parameter $\alpha$ behaves like the step size in gradient descent for the criterion weight $\grave{R}_k$. The first-order optimality criterion $r_k = 0$ must still be satisfied for the generic heuristic function $r_h({\rchi_{k{(-)}}})$ when the system is at the desired condition. Implicitly, this still optimizes some scalar objective; however, the exact objective need not be explicitly represented. This approach is also advantageous when implementing within the UKF framework, as the inverse Hessian is not required.

\section{Case Studies and Discussion} \label{section:case_studies_and_disc}

\subsection{Linear-Quadratic MPC with Observed Control} \label{subsec:msd_step_discussion}

A classic mass-spring-damper (MSD) system, with dynamics $m\ddot{x} + b\dot{x}+kx=F_u$ where $m=1$, $b=1$, and $k=2$ serves as a benchmark for evaluating the reactive, predictive, and anytime variants of observed control. The MPC objective is chosen to be the cost function \eqref{eqn:lqr_minimization} with weights
\begin{align}
    \label{MSD_objective_weights}
    Q_{p}= \begin{bmatrix}
        10 & 0 \\
        0 & 1
    \end{bmatrix} &&
    \begin{aligned}
        R_{p}  & = \begin{bmatrix}
            0.1
        \end{bmatrix}\\
        \tilde{R}_{p} & = \begin{bmatrix}
            0.1
        \end{bmatrix}
    \end{aligned} &&
    M_{p} = \begin{bmatrix}
        0\\
        0
    \end{bmatrix}.
\end{align}
Figure \ref{fig:msd_step_response_trajectories} compares Algorithms~\ref{alg:oc_plain}-\ref{alg:anytime-any-horizon MPC} performing step inputs for various time horizon lengths $N$, with a discrete LQR controller baseline; all $\Delta t = 0.1$. The discrete LQR control law is
\begin{equation}
    u_{lqr} = -K_{lqr}x_{ref} + u_{ss}
\end{equation}
where $u_{ss}=B^{\dagger}(\mathbb{I}-A)x_{ref}$ compensates for steady-state error caused by penalizing control effort. This term ensures convergence to steady state by solving $x_{ref} = A x_{ref} + B u_{ss}$
at the desired equilibrium $x_{ref}$, In observed control, the same behavior is achieved by making the measurement $H_k$
\begin{equation}
    H_k = \begin{bmatrix}
        \mathbb{I}_{n \times n} & \mathbb{O}_{n \times m}
    \end{bmatrix}
\end{equation}
which removes any control penalty.
Observed control's $u_{ss}$ is provided by the predictive time horizon, but remains bounded by the $\Delta u$ penalty, ensuring no integrator windup effects.

\begin{figure}[t!]
    \centering
    \includegraphics[width=0.5\textwidth]{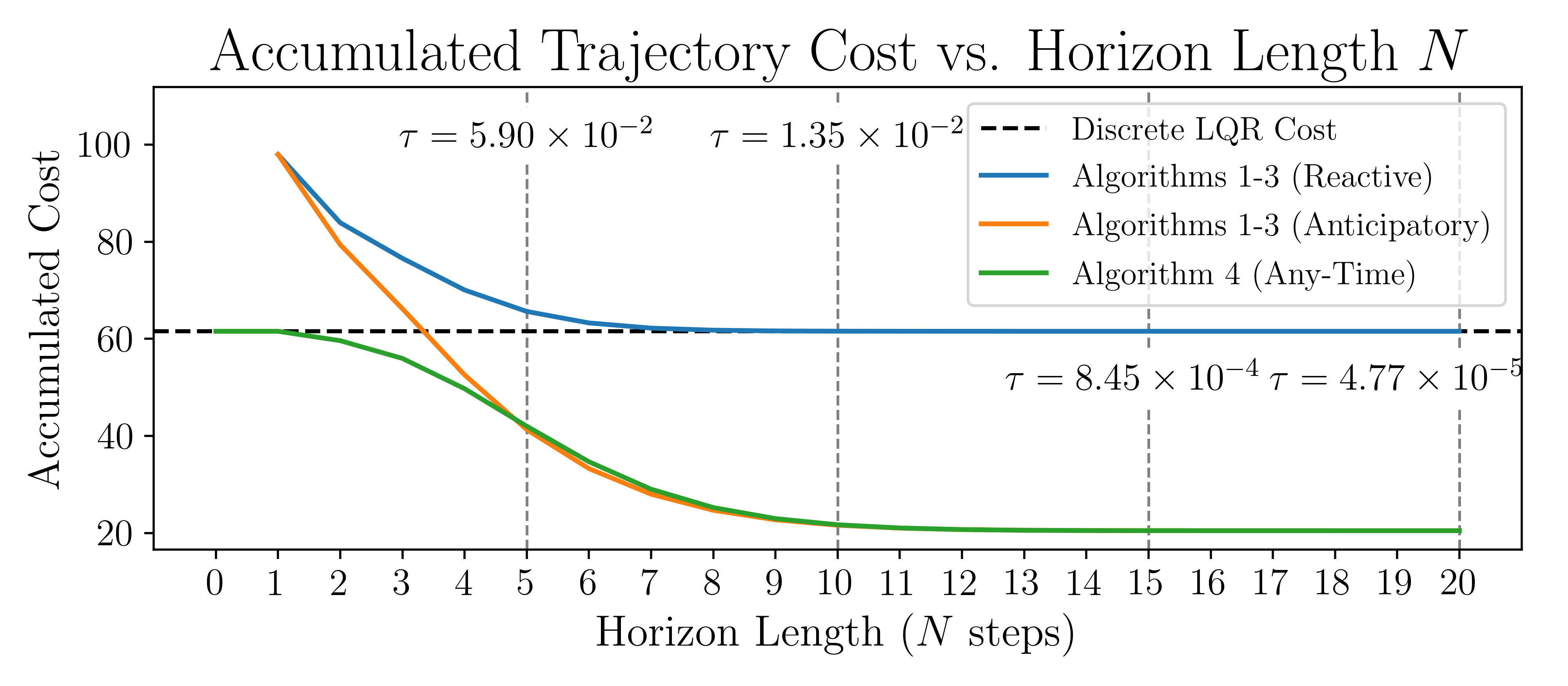}
    \caption{Accumulated trajectory costs for each observed control algorithm vs. horizon length $N$. All algorithms converge to the predictive optimal as $N$ increases with trajectory anticipation, with Algorithm \ref{alg:anytime-any-horizon MPC} starting at the reactive optimal with $N=0$. Conservatively, a $\tau$-termination threshold of $\delta_2=10^{-5}$ would ensure convergence. \vspace{-1em}}
    \label{fig:msd_horizon_costs}
\end{figure}

Figure \ref{fig:msd_horizon_reactive} shows reactive step responses of Algorithms \ref{alg:oc_plain}-\ref{alg:efficient_update} without trajectory anticipation, that is, each reference $(z_1, \cdots, z_N)$ is set to $z_0$, for various horizon lengths $N$. As expected, performance improves with larger $N$, ultimately converging to that of LQR. Figure \ref{fig:msd_horizon_anticipatory} shows predictive responses for Algorithms \ref{alg:oc_plain}-\ref{alg:efficient_update}, with larger $N$, the controllers optimally anticipate and center the step. Algorithm \ref{alg:anytime-any-horizon MPC}, the any-time any-horizon variant, shown in Figure \ref{fig:msd_horizon_anytime}, decomposes the horizon into purely reactive and anticipatory components, ensuring LQR's optimal reactive behavior for $N=0$, with predictive performance improving monotonically as $N$ increases.

Figure \ref{fig:msd_horizon_costs} presents the closed-loop trajectory cost accumulated over the entire trajectory given by
\begin{align}
    C_{\text{cl}} = \sum_{n=0}^{75} \left(\rchi_n^\top 
    \begin{bsmallmatrix}
        Q_{p} & M_{p}\\
        M_{p}^\top &  R_{p}
    \end{bsmallmatrix}
    \rchi_n + \|u_n-u_{n-1}\|^2_{\tilde{R}_{p}} \right)
\end{align}
with weights \eqref{MSD_objective_weights}
for each algorithm across various horizon lengths $N$. All controllers stabilize the system with no steady-state error, resulting in bounded accumulated costs. Algorithms \ref{alg:oc_plain}-\ref{alg:efficient_update} converge to the reactive optimum by $N=10$ whereas Algorithm \ref{alg:anytime-any-horizon MPC} begins at the reactive optimum with $N=0$, highlighting its advantage in short-horizon settings. By $N=12$, all algorithms converge to the predictive control optimum.

All algorithms demonstrate convergence to the predictive control optimum; however, the key contribution lies in Algorithm \ref{alg:anytime-any-horizon MPC} and its guaranteed baseline performance. It closes the gap for short time horizons, ensuring at worst LQR performance and stability, guaranteeing that there is always a reasonable stabilizing control ready while using any remaining compute time between updates for refinement.

\subsection{Non-Quadratic Costs and Measurement Modes} \label{subsec:meas_modes_discussion}

Quadratic objectives have well-defined minima; however, they limit the types of behaviors that can be encoded into a MPC controller. Table \ref{table:comp_agg_meas_summary} summarizes the construction of measurement residuals, sensitivities, and covariances for the non-quadratic case, which are equivalent to the duality-based formulation when the objective function is itself quadratic.

To illustrate the behavior of different measurement modes under non-quadratic objectives, consider the 2D linear system
\begin{equation}
    \label{eqn_lds_system_dynamics}
    m \begin{bmatrix}
        \ddot{x} & \ddot{y}
    \end{bmatrix}^\top + b \begin{bmatrix}
        \dot{x} & \dot{y}
    \end{bmatrix}^\top = \begin{bmatrix}
        F_x &
        F_y
    \end{bmatrix}^\top
\end{equation}
with $m=1$ and $b=1$. The controller is tasked with tracking a full-state reference trajectory $\rchi_{ref}$ while also avoiding reference-intersecting obstacles. The augmented state $\rchi$ is
\begin{align}
    \rchi = \begin{bmatrix}
        x & y & v_x & v_y & F_x & F_y
    \end{bmatrix}^\top
\end{align}
which is penalized by $W_\rchi$. To impose obstacle avoidance, the obstacle potential field cost $\mathcal{C}_\text{obs}(d)$ is introduced, given by
\begin{align*}
    \mathcal{C}_\text{obs}(d) =
    \begin{cases}
        \mathcal{V}(d_t) + \sum_{n=1}^{2} \frac{\partial^n\mathcal{V}}{\partial d^n} (d-d_t)^n & \text{if } 0 \leq d < d_t\\
        \mathcal{V}(d)& \text{if } d_t \le d < d_z \\
        0 & \text{if } d_z \le d
    \end{cases}
\end{align*}
where $d_t=0.1$ and the underlying potential field $\mathcal{V}(d)$ is
\begin{align*}
    \mathcal{V}(d) &= (1/100)\left( \tan(i(d))/i(d) - i(d)^2/3-1\right)\\
    \notag
    i(d) &= -(\pi/2)\left( 1 - (d-d_c)/d_z \right)
\end{align*}
where $d(\rchi)$ is the closest obstacle distance and $d_z = 0.5$ specifies the zero-cost distance. The function $\mathcal{C}_\text{obs}(d)$ is constructed such that it has continuous first and second-order derivatives. The generic objective function stage contribution $J_k$ of \eqref{eqn:disc_time_opt_ctrl_formulation:obj_fn} is
\begin{align}
    \label{eqn_lds_obs_scalar_objective}
    J_k = (1/2)\| \rchi_{k, ref} - \rchi_{k} \|_{W_\rchi}^2 + W_\text{o} \mathcal{C}_\text{obs}(d)
\end{align}
where $W_\rchi = \mathrm{diag}([100, 100, 1, 1, 0.1, 0.1])$ and $W_\text{o}=20000$; note the $\Delta u$ term in \eqref{eqn:disc_time_opt_ctrl_formulation:obj_fn}, weighted $\tilde{R}_{p}=0.1$.
The potential field term $\mathcal{C}_\text{obs}(d)$ can produce an indefinite Hessian; to avoid ill-posedness from its contributions in the hessian of \eqref{eqn_lds_obs_scalar_objective}, the complementary projector $T$, constructed from the eigen-decomposition $\frac{\partial^2 \mathcal{C}_\text{obs}}{\partial \rchi^2} = V \Lambda V^\top$, is used to create the semi-definite $C_{x}$ and $C_{xx}$ for computing $J_{x}$ and $J_{xx}$ of \eqref{eqn_lds_obs_scalar_objective}
\begin{align*}
    \mathcal{C}_{x} \gets T \frac{\partial \mathcal{C}_\text{obs}}{\partial \rchi} &&
    \mathcal{C}_{xx} \gets T \frac{\partial^2 \mathcal{C}_\text{obs}}{\partial \rchi^2} T^\top
    &&
    T = \mathbb{I} - \sum_{\lambda_i < 0} v_i v_i^\top.
\end{align*}

\begin{figure}[t!]
    \vspace{-1.5em}
    \centering
    \includegraphics[width=0.5\textwidth]{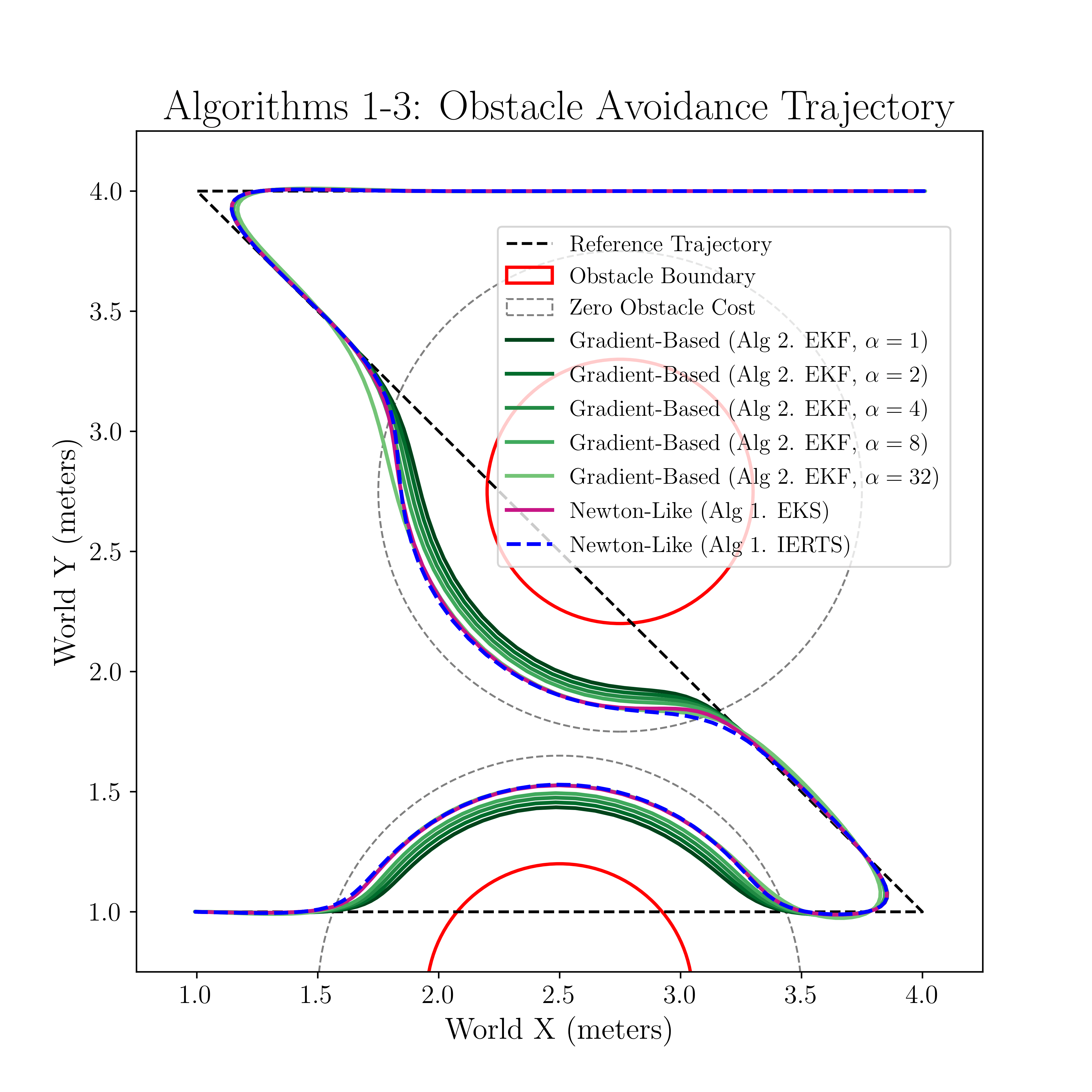}
    \caption{Observed control performing obstacle avoidance on a collision-containing trajectory with a nonquadratic objective \eqref{eqn_lds_obs_scalar_objective} and system \eqref{eqn_lds_system_dynamics}. The Newton-like mode with IERTS finds the local optimal while with EKS approximates the locally optimal trajectory in one step. The gradient-based mode with various $\alpha$ values demonstrates the tunability of \textit{direction-to-better} objectives. \vspace{-1.5em}}
    \label{fig:lds_osbtacle_trajectories}
\end{figure}

Figure \ref{fig:lds_osbtacle_trajectories} illustrates the behaviors resulting from the Newton-like and gradient-based $(\alpha = 1, 2, 4, 8, 32)$ measurement modes for various nonlinear backends. Using \eqref{numerically_safe_newton_step}, under the Newton-like measurement mode, the posterior-prior difference $\widetilde{\rchi}_k = \rchi_{k(+)} - \rchi_{k(-)}$ at each step can be written as
\begin{align}
    \label{eqn_lds_obs_scalar_map}
    \widetilde{\rchi}_k = - P_{k(-)}\left(P_{k(-)} + J_{xx}^{-1} \right)^{-1} r_k.
\end{align}
Algorithm \ref{alg:oc_foward_only} highlights how at each step, the posterior-prior differences map through smoother gains to become control refinements and \eqref{eqn_lds_obs_scalar_map} highlights how the ratio between the covariance $P_{k(-)}$ and the stage cost dual $J_{xx}^{-1}$ weight the residual direction $r_k$. Likewise \eqref{eqn_newton_measurement_information_form} reveals Algorithm \ref{alg:oc_foward_only} with the Newton-like measurement mode computes the control from accumulated information-regularized Newton steps.

\begin{figure*}[t!]
    \centering
    \begin{minipage}[t]{0.65\textwidth}
        \centering
        \begin{minipage}[t]{0.485\linewidth}
            \centering
            \includegraphics[width=\linewidth]{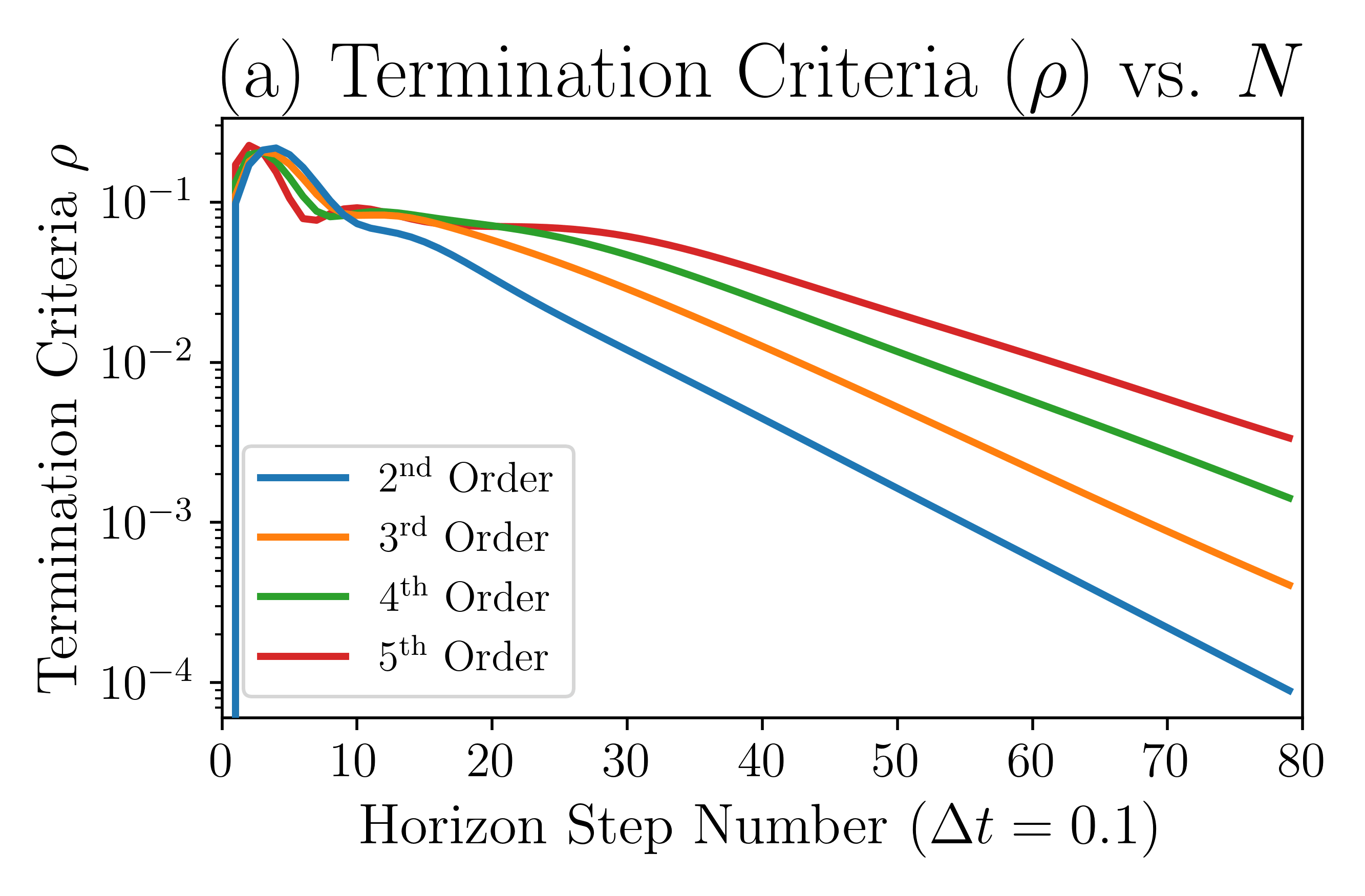}
        \end{minipage}%
        \begin{minipage}[t]{0.485\linewidth}
            \centering
            \includegraphics[width=\linewidth]{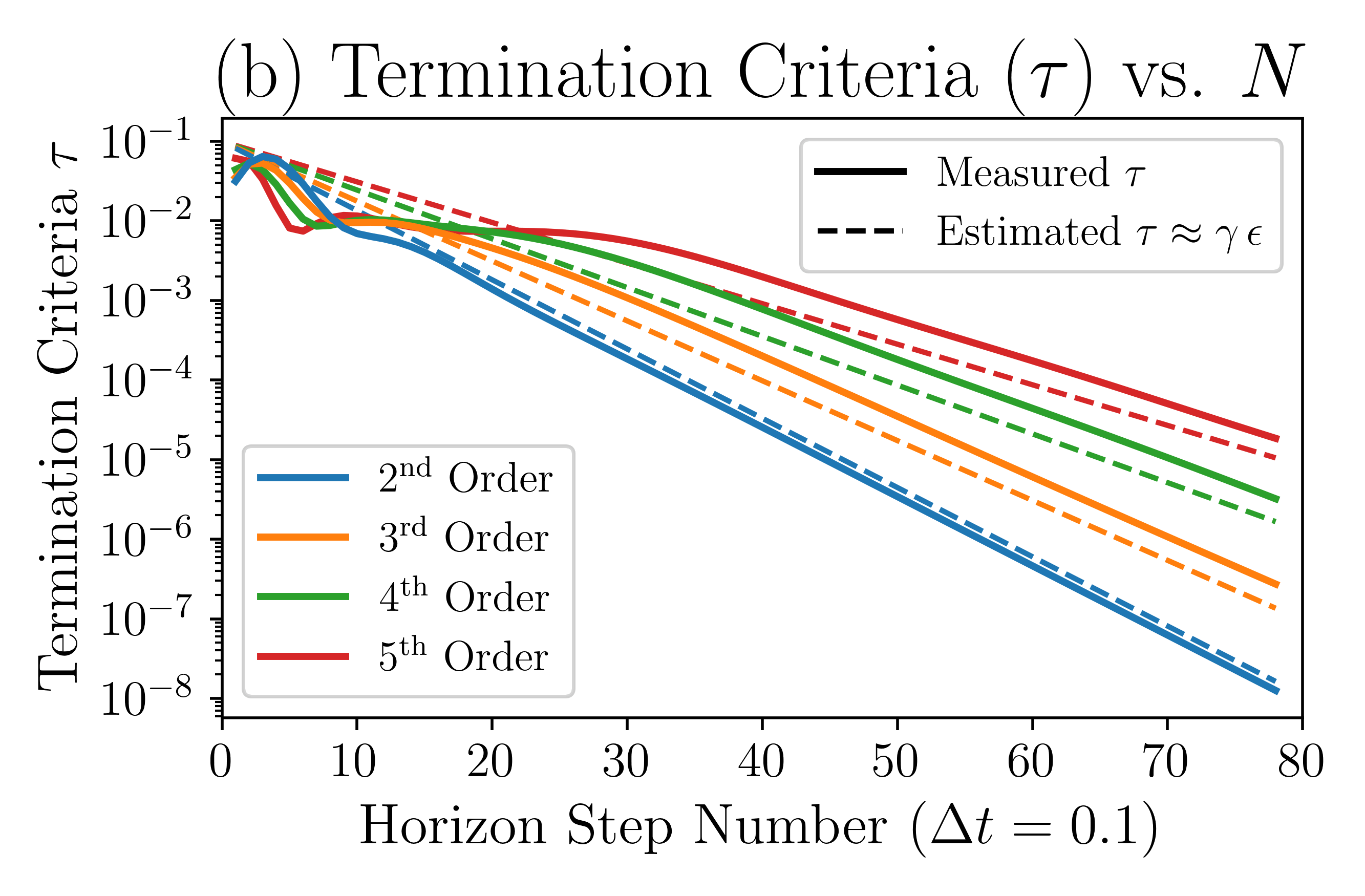}
        \end{minipage}
        \captionof{figure}{The $\rho$ and $\tau$ with $(\gamma=0.1)$ early termination criteria vs horizon step number $N$ for various order linear systems. Notably, the accumulator gain criteria $\rho$ peaks when the system first becomes fully controllable. The $\tau$ estimate is found by rearranging \eqref{eqn_required_N_estimate}. \vspace{-1em}}
        \label{fig:termination_criteria_side_by_side}
    \end{minipage}
    \hfill
    \begin{minipage}[t]{0.32\textwidth}
        \centering
        \includegraphics[width=\linewidth]{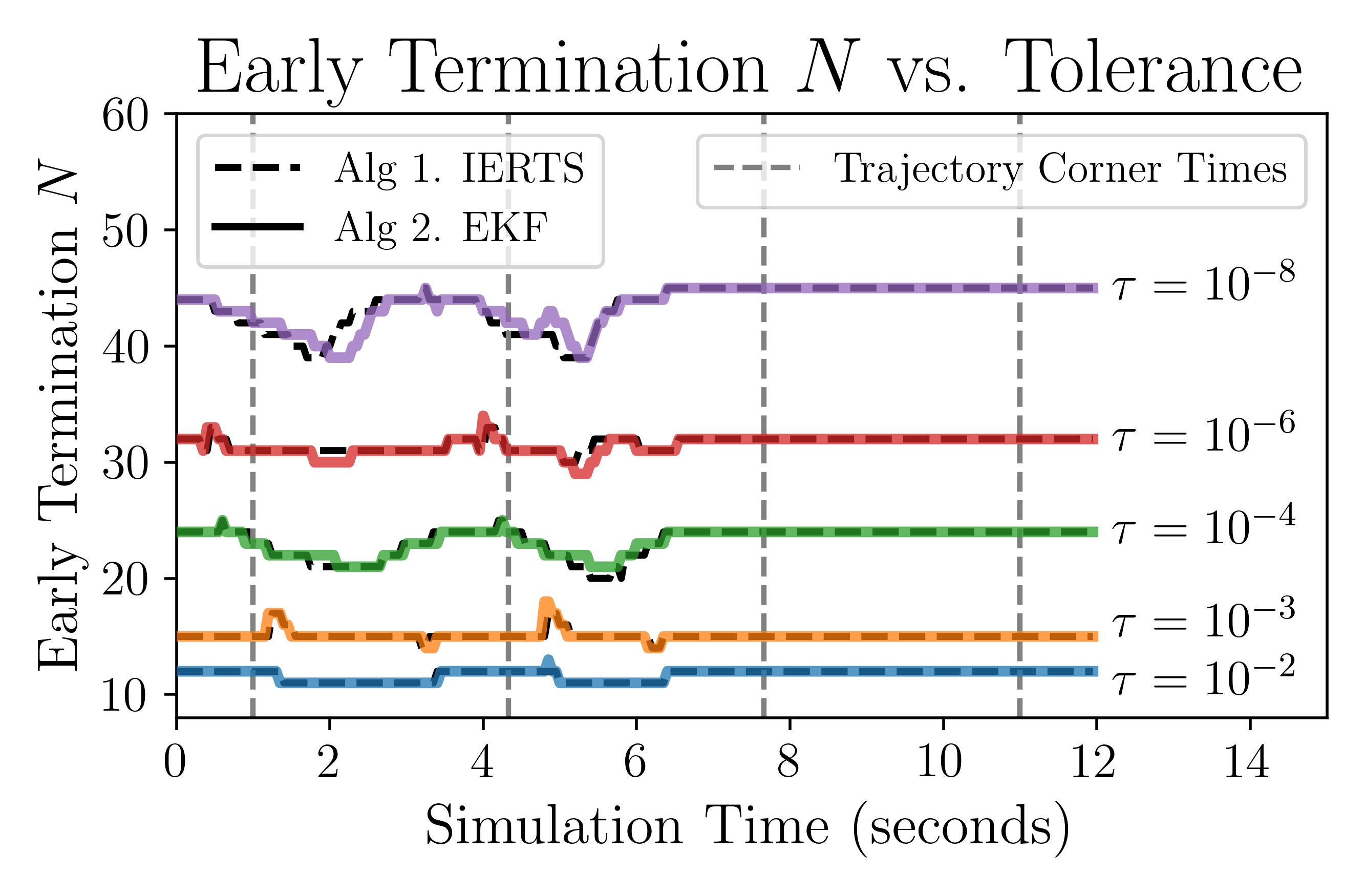}
        \captionof{figure}{Terminated horizon lengths for the obstacle avoidance scenario (Newton-like measurement mode) for various $\tau$ thresholds.}
        \label{fig:lds_termination_lengths}
    \end{minipage}
    \vspace{-1em}
\end{figure*}

However, when tuning an MPC controller, it is often easier to describe a \textit{direction to better} $(\nicefrac{\partial J}{\partial \rchi_k})$ and \textit{how much to care} weight $\alpha$ for which gradient-based measurements provide for a more intuitive interface. In the gradient framework, the objective is specified as a residual function $r_h$, and weight $\grave{R}$
\begin{align*}
    r_h(\rchi_k)= \begin{bmatrix}
        \rchi_{k, ref} - \rchi_{k}  \\
        \frac{\partial \mathcal{C}}{\partial \rchi}
    \end{bmatrix} &&
    \grave{R} = \alpha \begin{bmatrix}
        W_\rchi^{-1} & 0 \\
        0 & W_o^{-1}
    \end{bmatrix}
\end{align*}
and the posterior-prior difference is given by
\begin{align*}
    \widetilde{\rchi}_k
    = P_{k(-)}  \frac{\partial r_h}{\partial \rchi_{k}}^\top \Biggl[
        \frac{\partial r_h}{\partial \rchi_{k}} P_{k(-)} \frac{\partial r_h}{\partial \rchi_{k}}^\top  + \grave{R}
    \Biggr]^{-1} r_h(\rchi_k)
\end{align*}
where $\alpha$ can be intuited as a gradient descent step size in the ratio of how $r_h(\rchi_k)$ maps to control refinements.
The gradient-based method supports a broader class of objectives and behaviors---including those without an analytical $J(\rchi)$, those impractical to directly evaluate, or those incorporating learned objectives or heuristics---and offers greater tunability.

\subsection{Early Termination}
\label{subsec:early_term_discussion}

The key feature of Algorithms \ref{alg:oc_foward_only}-\ref{alg:anytime-any-horizon MPC} is their single forward-only structure, iterating through horizon steps sequentially into the predictive horizon. This structure enables the controllers to precisely quantify the additional control refinement on $u_0$ resulting from extending the predictive horizon length from $N$ to $N+1$, via the termination metrics $\rho$ and $\tau$ introduced in Section \ref{subsec:oc_early_term_crit}. These metrics can be precomputed for linear systems or evaluated on-the-fly using repeated linearizations and the filter's covariances for nonlinear systems enabling adaptive horizon length selection and early termination of the MPC optimization without sacrificing control performance.

Figure \ref{fig:termination_criteria_side_by_side} illustrates both termination criteria as a function of the horizon length $N$ for several representative linear systems. Each $n$-th order system is constructed with the first $n$ poles chosen from the set $\mathcal{P}=(-0.1, -0.2, -0.3, -0.4, -0.5)$ with the system matrices $A$ and $B$ defined as
\begin{align}
    \label{termination_criteria_system}
    A = \begin{bmatrix}
            \begin{matrix}                
                \vec{0}_{n-1}  & & \mathbb{I}_{n-1}
            \end{matrix} \\
            \begin{bmatrix}
                -a_0 & -a_1 & \cdots & -a_{n-1}
            \end{bmatrix}
    \end{bmatrix} &&
    B = \begin{bmatrix}
        \vec{0}_{n-1}\\
        1
    \end{bmatrix}
\end{align}
where $\{ a_0, a_1, \cdots, a_{n-1} \}$ are the coefficients of the characteristic polynomial defined by the chosen poles. This construction ensures each system's controllability index matches its order.

As shown in Figure \ref{fig:termination_criteria_side_by_side}a, the accumulator gain criteria $\rho$ decreases as additional future-time horizon steps have diminishing influence on $u_0$. Small periodic dips appear after controllability is established and settle thereafter, corresponding to the transient convergence behavior being reflected in the accumulators $\Gamma_k$ and $G_k$ used in finding $\rho$.
Similarly, Figure \ref{fig:termination_criteria_side_by_side}b shows the uncertainty metric $\tau$ and the horizon estimation heuristic \eqref{eqn_required_N_estimate} as a function of $N$. For this system, the convergence ratio for $P_{\infty[s]}/K_\infty$ was empirically determined to be $\gamma \approx 0.1$. Here, the transient $\tau$-dynamics are also captured early, and $\tau$ exhibits geometric (exponential) convergence consistent with the heuristic \eqref{eqn_required_N_estimate}. This confirms that the controllers rapidly achieve the minimum horizon length required to sufficiently approximate the infinite-horizon $u_0$ under a stationary reference trajectory. 

For the obstacle avoidance problem with the system \eqref{eqn_lds_system_dynamics}, Figure \ref{fig:lds_termination_lengths} shows the adaptive horizon lengths selected throughout the entire simulation as the system is controlled from the lower left corner to the upper right corner. The Newton-like measurement is used, and the $\rho$ metric ignored to isolate the effects of $\tau$ as references are known smooth. As the cost gradient becomes better defined, particularly near obstacles, the number of horizon steps required to reach a given control confidence decreases. Predictably, tightening the tolerance increases the average number of required horizon steps, but the curve retains its general shape. The approximately optimal controllers produced by EKFs determines very similar horizon lengths compared to the converged IERTS-based controllers.

Even for linear systems, selecting an appropriate horizon is difficult; for nonlinear systems, nonquadratic objectives, or nonstationary trajectories, the difficulty compounds. 
Corollary \ref{cor: delta epsilon} and \ref{cor: delta p0 -> 0} emphasize the importance of jointly considering both convergence to an infinite-time approximation and impact of reference trajectory changes, while Corollary \ref{cor:smoothness_of_trajectory} connects horizon selection to reference trajectory smoothness. By quantifying, at runtime, the incremental refinement in $u_0$ from each additional horizon step, observed control enables the dynamic determination of minimum-length prediction horizons. Although, as shown in Figure \ref{fig:lds_termination_lengths}, approximately optimal controllers generally closely match fully-converged controllers in horizon length selection, the $\rho$ and $\tau$ termination metrics are only exact for fully-converged horizons.

\subsection{Non-Linear Cart Pole Swing Up} \label{subsec:non_linear_system}

\begin{figure}[t!]
    \centering
    \includegraphics[width=.5\textwidth]{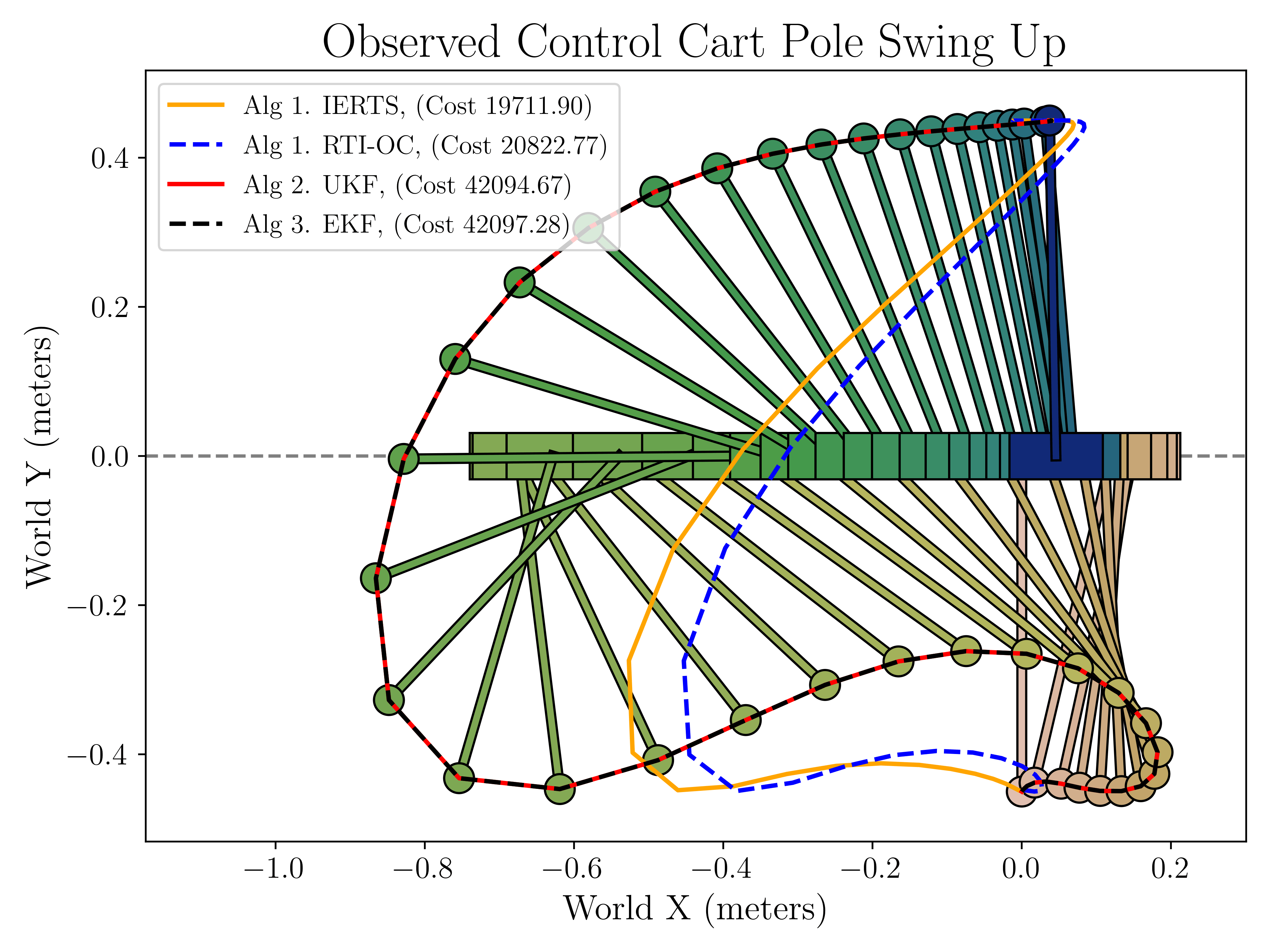}
    \caption{Cart pole swing-ups performed with different estimation backends (Newton-like measurement, $\Delta t=0.05$, $N=250$). The EKF- and UFK-based approximate controllers produce effectively indistinguishable trajectories. Notably, RTI-OC computes a very similar solution to IERTS with a fraction of the compute time.\vspace{-1.5em}}    \label{fig:swingup}
\end{figure}

Figure \ref{fig:swingup} shows observed control performing a cart-pole swing up --- a widely used benchmark for nonlinear predictive control. In this experiment, force control was applied to the cart's horizontal axis, and the pole is modeled as a weightless rigid link with all the mass concentrated at the endpoint. The augmented state $\rchi$ for the system is defined as
\begin{equation}
    \label{eqn_cartpole_state_def}
    \rchi = \begin{bmatrix}
        x & v & \theta & \omega & F_x
    \end{bmatrix}^\top
\end{equation}
where $x$ is the cart's horizontal position, $v=\dot{x}$ is the cart's horizontal velocity, $\theta$ is the pole angle relative to the vertical (counter-clockwise positive), $\omega=\dot{\theta}$ is the angular velocity of the pole, and $F_x$ is the control input force applied to the cart. The nonlinear dynamics of the system are given by
\begin{subequations}
    \begin{alignat}{2}
        \ddot{x} &~=\,&& \left[
            l ( m_p \sin^2\theta + m_c )
        \right]^{-1} \left[
            l F_x - d_v l \dot{x} \right.\\
        & && \qquad \left. - m_p l^2 \dot{\theta}^2 \sin\theta - \left( d_\omega \dot{\theta} - g m_p l \sin\theta \right)\cos\theta
        \right] \nonumber \\
        \ddot{\theta} &~=\,&& \left[
            m_p l^2 ( m_p \sin^2\theta + m_c )
        \right]^{-1}  \\
        & && \qquad \times \left[
            \left( m_p l F_x - d_v m_p l \dot{x} \right) \cos\theta - d_\omega ( m_c + m_p ) \dot{\theta} \right. \nonumber \\
        & && \qquad \left. + \left( 
                - m_p^2 l^2 \dot{\theta}^2 \cos\theta
                + g m_p^2 l + g m_c m_p l 
            \right) \sin\theta 
        \right] \nonumber
    \end{alignat}
\end{subequations}
where the parameters are: cart mass $m_c=0.25$kg, pole mass $m_p=0.2$kg, pole length $l=0.45$m, linear damping $d_v=0.05$Ns/m, and angular damping $d_\omega=0.015$Nms/rad. The objective function used is quadratic and is given by
\begin{align}
    \label{eqn_cart_pole_obj}
    J_{cp} = \sum_{k=0}^N  \left( \| \rchi_k - \rchi_{ref,k}\|^2_{W_{cp}} + (\Delta F_{x,k})^2 W_{f} \right)
\end{align}
where the weight $W_{cp} = \mathrm{diag}([1000, 1, 300, 25])$ penalizes the state \eqref{eqn_cartpole_state_def} and the weight $W_{f}=20$ penalizes changes in control effort. The reference $\rchi_{ref,k}$ is defined as the step
\begin{align}
    \rchi_{ref, k} = [0 ~~ 0 ~~ \pi ~~ 0 ~~ 0 ]^\top & \text{if } t < 1; \mathbb{O}_{5\times 1} \text{ otherwise}.
\end{align}

Figure \ref{fig:swingup} shows observed control performing a cart pole swing up with different backends. With EKF and UKF backends, the results are consistent with the general understanding that for moderate nonlinearities, both EKF and UKF provide reasonably accurate covariance propagation, and thus in observed control, result in nearly equal controls. The IERTS backend provides locally optimal controls, but requires iterating. Performing RTI-OC, which hot-starts and performs only one iteration per control update, results in a similar solution to IERTS but with far less compute. Initializing RTI-OC with a single converged IERTS update would eliminate the early swing up deviation, yielding IERTS-like performance with the same one-smoother-pass cost after the first update.

\begin{figure}[t!]
    \centering
    \includegraphics[width=.5\textwidth]{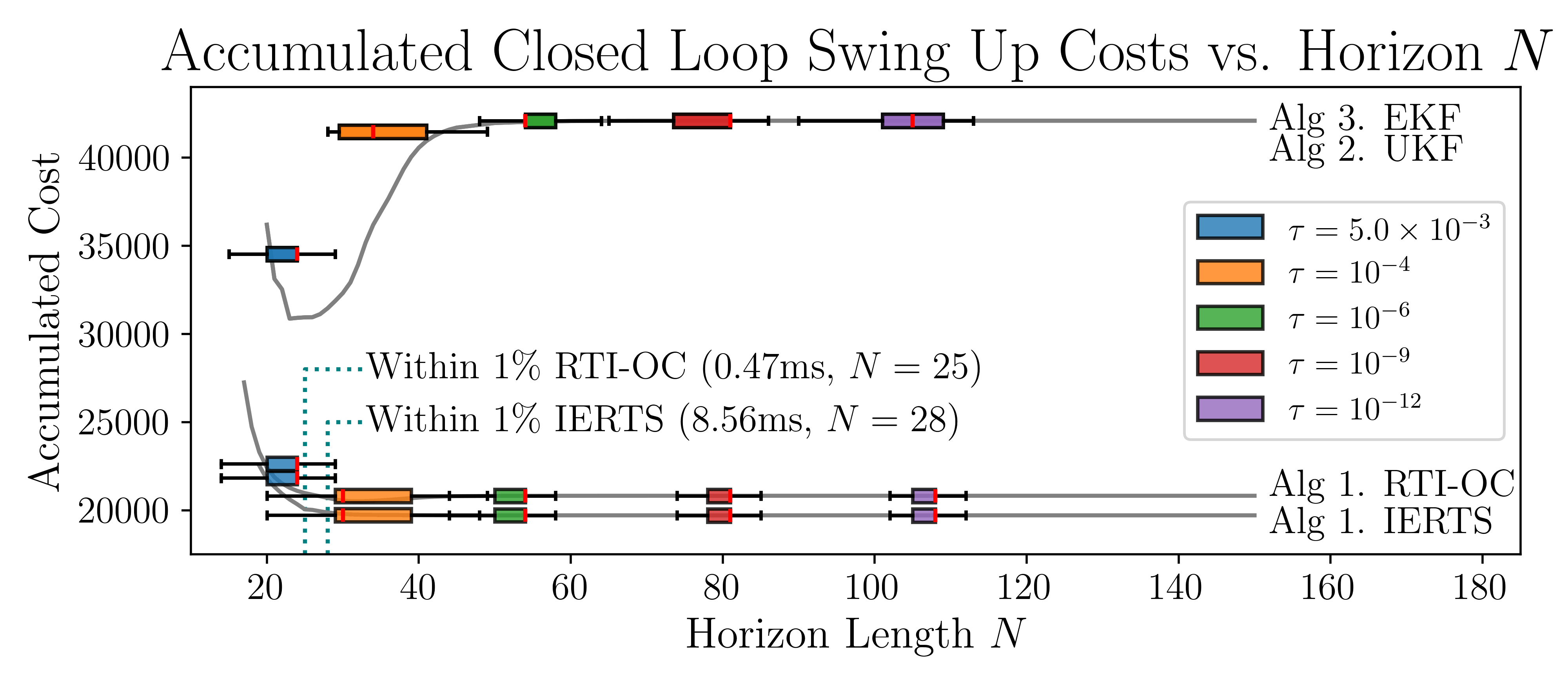}
    \caption{Controller objective function \eqref{eqn_cart_pole_obj} performance vs horizon lengths for different backends. IERTS finds the local optima, to which RTI-OC comes close, while the approximate EKF- and UKF-based controllers are best-performing at $N=23$. Vertical lines indicate when the controller is within 1\% of its $N=500$ result and the box plots outline the adaptive horizon lengths selected when early termination is enabled for the specified $\tau$ tolerance. \vspace{-1em}}
    \label{fig:swingup_cost_vs_N}
\end{figure}

Figure \ref{fig:swingup_cost_vs_N} shows the accumulated objective function performance of observed control for various backends and horizon lengths. The locally optimal IERTS backend exhibits the expected behavior, i.e., decreasing in cost as $N$ increases. By $N = 28$, IERTS is within 1\% of the $N=500$ result; a performance level that a $\tau$-termination tolerance of $\tau\le 10^{-4}$ would achieve. Similarly, RTI-OC experiences decreasing accumulated costs with increasing $N$, though it is slightly suboptimal due to its incomplete convergence each update. Notably, RTI-OC is 20x faster than IERTS and achieves sub-millisecond converged solutions. When using EKF and UKF, the costs first decrease to $N=23$, then increase as the filters increasingly low-pass-filter future stabilizing information resulting in conservative, yet suboptimal, control actions.

\begin{figure}[t!]
    \centering
    \begin{subfigure}[b]{0.5\textwidth}
        \centering
        \includegraphics[width=\textwidth]{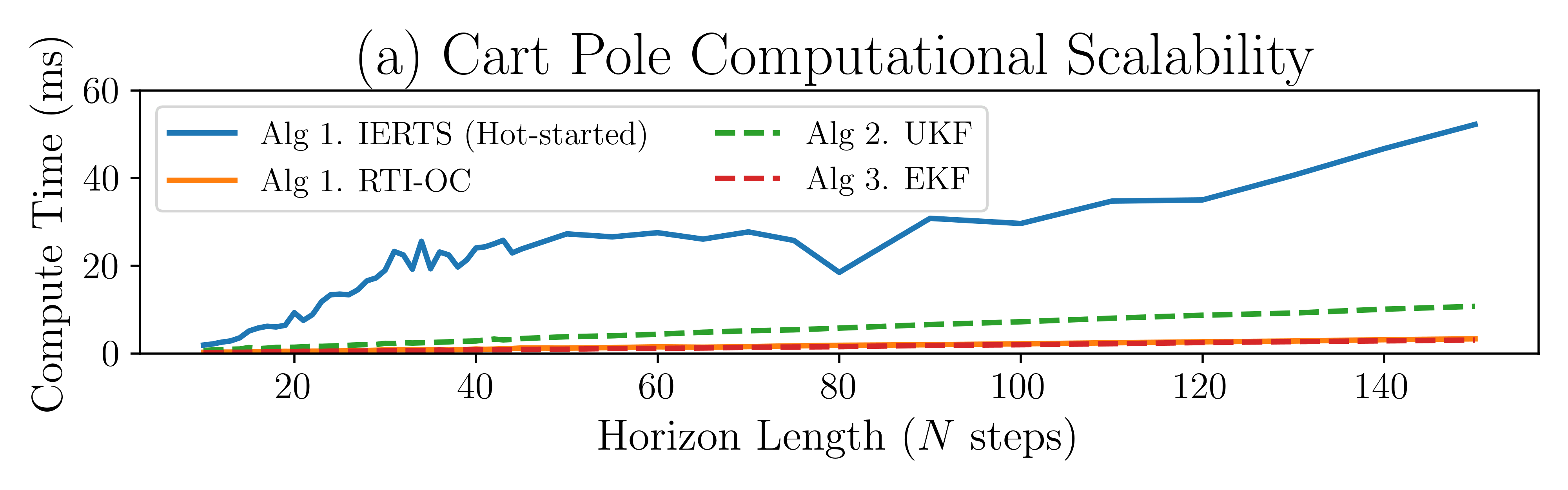}
        \phantomsubcaption%
        \vspace{-1.5em}
        \label{fig:swingup_time_vs_N}
    \end{subfigure}
    \begin{subfigure}[b]{0.5\textwidth}
        \centering
        \includegraphics[width=\textwidth]{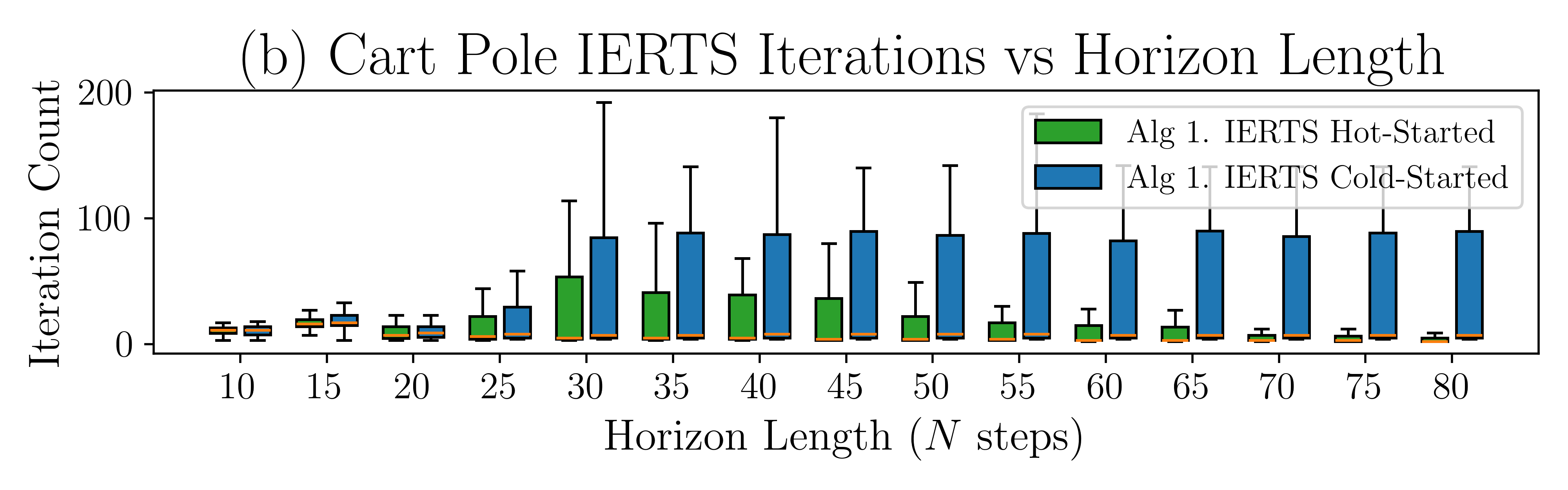}
        \phantomsubcaption%
        \vspace{-1.5em}
        \label{fig:swingup_iterations_vs_N}
    \end{subfigure}
    \caption{
    Average computational time per control update (a) for each backend and iteration count (b) for the IERTS backend. Observed control is linearly scalable with respect to the horizon length.
    \vspace{-1.5em}}
    \label{fig:cp_cost_and_scalability_figure}
\end{figure}

Computational performance for the cart-pole swing up task is benchmarked in Figure \ref{fig:swingup_time_vs_N}. All benchmarks were written in C++20, used Eigen 3.4.0, compiled with GCC 13.3, and ran on an Intel i7-12700K CPU (single threaded). Figure \ref{fig:swingup_iterations_vs_N} shows how the IERTS iteration count varies with prediction horizon length $N$. For cold-started IERTS, the iterations count plateaus as $N$ increases due to the strong temporal band-structure of the MPC problem. For hot-started IERTS, the iteration count instead decreases after $N\approx35$ as the augmented states become increasingly converged as they shift through the horizon. Computational scalability is shown up to $N=150$, though this greatly exceeds the horizon length necessary to be within 1\% of optimality. As such, a threshold of between $10^{-6}$ and $10^{-8}$ is recommended for $\tau$. Together, these results show that observed control is a highly effective and computationally scalable approach to nonlinear unconstrained MPC. Extending the framework to accommodate general general state and input constraints remains an important direction for future work.

\section{Conclusion} \label{section:conclusions}

Observed control leverages duality between state estimation smoothers and model predictive control to achieve linear time-horizon scalability and adaptive horizons. Any filter that is compatible with RTS smoothers and supports zero process noise can be an observed control backend. The presented algorithms, based on Kalman-type filters, inherit their properties---providing optimal controls for linear systems with quadratic penalties and high-quality or optimal controls for nonlinear systems or non-quadratic penalties (depending on the underlying estimator), and demonstrate millisecond-level single-threaded compute times. Thus, observed control is an attractive MPC algorithm that enables dual-use of existing state-estimation codes for MPC computation that is appropriate for high-rate, resource-constrained, nonlinear predictive control. 

\section{References}

\bibliographystyle{unsrt}
\bibliography{references}

\end{document}